\documentclass[12p]{article}
\RequirePackage{amsthm,amsmath,amsfonts,amssymb}
\usepackage{xcolor}
\usepackage{setspace}
\usepackage[left=2.5cm,right=2.5cm,top=2.5cm,bottom=2cm]{geometry}
\usepackage{graphicx}
\usepackage[authoryear, round]{natbib}
\usepackage{tikz}
\usetikzlibrary{angles,quotes}
\usepackage{pgfplots}

\newtheorem{theorem}{Theorem}
\newtheorem{lemma}{Lemma}

\newcommand{\mb}{\mathbf}
\newcommand{\bX}{\mathbf{X}}
\newcommand{\bx}{\mathbf{x}}

\newcommand{\by}{\mathbf{y}}
\newcommand{\bP}{\mathbf{P}}

\newcommand{\bS}{\mathbf{S}}

\newcommand{\bu}{\mathbf{u}}

\newcommand{\cF}{\mathcal{F}}
\newcommand{\rmi}{\mathrm{i}}
\newcommand{\eps}{\varepsilon}

\newcommand{\E}{\mathbb{E}}  
\newcommand{\bal}[1]{\begin{align*}#1\end{align*}}

\begin{document}

\title{The smallest singular value of large random rectangular Toeplitz and circulant matrices}
\author{Alexei Onatski and Vladislav Kargin}

\maketitle
\begin{abstract}
    Let $x_i$, $i\in\mathbb{Z}$ be a sequence of i.i.d.~standard normal random variables. Consider rectangular Toeplitz $\bX=\left(x_{j-i}\right)_{1\leq i\leq p,1\leq j\leq n}$ and circulant $\bX=\left(x_{(j-i)\mod n}\right)_{1\leq i\leq p,1\leq j\leq n}$ matrices. Let $p,n\rightarrow\infty$ so that $p/n\rightarrow c\in(0,1]$. We prove that the smallest eigenvalue of $\frac{1}{n}\bX\bX^\top$ converges to zero in probability and in expectation. We establish a lower bound on the rate of this convergence. The lower bound is faster than any poly-log but slower than any polynomial rate. For the ``rectangular circulant'' matrices, we also establish a polynomial upper bound on the convergence rate, which is a simple explicit function of $c$.  
\end{abstract}

\section{Introduction} 

Let $x_i$, $i\in\mathbb{Z}$ be a sequence of i.i.d.~standard normal random variables. Consider a $p\times n$ random Toeplitz matrix with $n\geq p$
\begin{equation}
\label{toeplitz matrix}
\bX=\begin{pmatrix}
    x_{0}&x_{1}&x_{2}&&\cdots&x_{n-2}&x_{n-1}\\
    x_{-1}&x_{0}&x_{1}&\ddots&&x_{n-3}&x_{n-2}\\
    x_{-2}&x_{-1}&x_{0}&\ddots&&x_{n-4}&x_{n-3}\\
    \vdots&\ddots&\ddots&\ddots&&\vdots&\vdots\\
    x_{1-p}&x_{2-p}&\dots&&&x_{n-1-p}&x_{n-p}
\end{pmatrix}.
\end{equation}
The $(i,j)$-th entry of $\bS:=\frac{1}{n}\bX\bX^\top$ can be viewed as a sample autocovariance at lag $|i-j|$ based on the sample $x_{1-p},\dots,x_{n-1}$. If only the observations with non-negative indices are available, one may replace $x_{1-p},x_{2-p},\dots,x_{-1}$ in the above definition of $\bX$ by $x_{n+1-p},x_{n+2-p},\dots,x_{n-1}$, respectively, to obtain a \textcolor{black}{$p\times n$} ``rectangular circulant'' matrix
\begin{equation}
\label{circulant matrix}
\bX=\begin{pmatrix}
    x_{0}&x_{1}&x_{2}&&\cdots&x_{n-2}&x_{n-1}\\
    x_{n-1}&x_{0}&x_{1}&\ddots&&x_{n-3}&x_{n-2}\\
    x_{n-2}&x_{n-1}&x_{0}&\ddots&&x_{n-4}&x_{n-3}\\
    \vdots&\ddots&\ddots&\ddots&&\vdots&\vdots\\
    x_{n+1-p}&x_{n+2-p}&\dots&&&x_{n-1-p}&x_{n-p}
\end{pmatrix}.
\end{equation}
Then the entries of $\bS$ become the so called circular sample autocovariances  \cite[ sec.~6.5.2]{anderson1971}. 

The sample autocovariance and circular autocovariance matrices and their inverses are fundamental for time series analysis. The norm of $\bS^{-1}$, and hence the smallest eigenvalue of $\bS$, plays an important role in the studies of autoregressive spectral density estimates and optimal prediction (e.g.~\cite{berk73}, \cite{bhansali78,bhansali96}, \cite{shibata80},  \cite{ing03,ing05}). 

Random Toeplitz and circulant matrices also play an important role in the randomized preprocessing of linear systems of equations (\cite{pan13}) and in compressed sensing (\cite{rauhut09}). In these fields, the smallest eigenvalue of $\bS$ is particularly important because it influences the condition number, which determines the accuracy of various numerical procedures.



This paper shows that the smallest eigenvalue of $\bS$, for both the Toeplitz and circulant cases, converges to zero in probability and in expectation when $p,n\rightarrow\infty$ so that $p/n\rightarrow c\in(0,1]$. 
The convergence to zero happens no matter how small $c>0$ is, which stays in sharp contrast to the standard sample covariance case, where all elements of $\bX$ are i.i.d.  
As is well known (see \cite{Silverstein85} for the first such result), in the latter case, the smallest eigenvalue of $\bS$ almost surely converges to $(1-\sqrt{c})^2$ provided that $c\leq 1$. As $c$ approaches zero, this limit approaches $1$, the population variance of $x_i$. Note that, similar to the case of sample covariance matrices, any two rows of a random Toeplitz or circulant matrix $\bX$ are nearly orthogonal. Therefore, some readers may find the convergence of the smallest eigenvalue to zero \textit{for any} $c\in (0,1]$ surprising, as we do.

Moreover, we show that the convergence to zero is faster than any power of the logarithm of $n$ (poly-log rate). For the ``rectangular circulant'' matrix, we show that, conversely, the convergence is slower than any power of $n$ (polynomial rate) if $c$ is sufficiently small. This slow rate makes it difficult to verify, or even guess, the convergence using Monte Carlo (MC) experiments.

\begin{figure}
\includegraphics[]{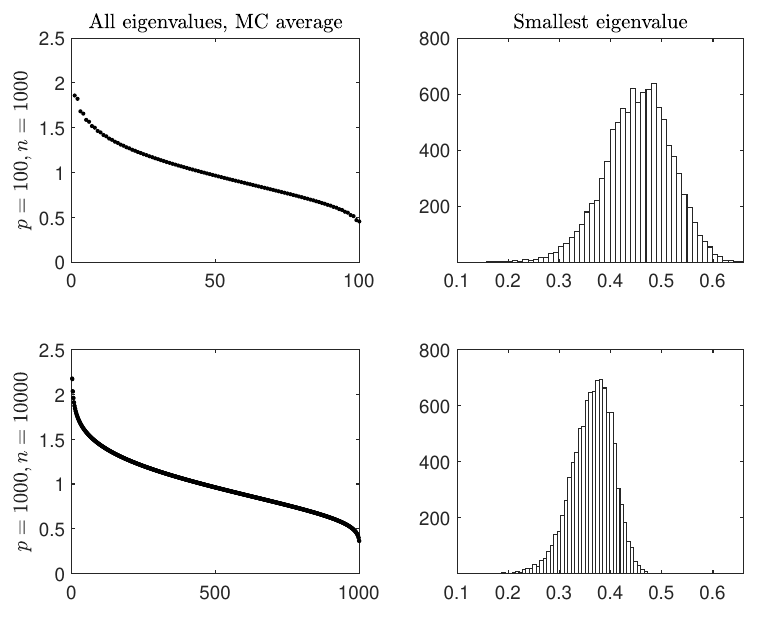}
\centering
    \label{fig:MCillustration}
    \caption{Monte Carlo results for random circulant matrix \eqref{circulant matrix}. The number of MC replications is 10,000. Top panel: $p=100, n=1,000$. Bottom panel: $p=1,000,n=10,000$. Left panel: MC average of all the eigenvalues of $\bS$. Right panel: histogram of MC realizations of the smallest eigenvalue. }
\end{figure}

To illustrate this difficulty, Figure \ref{fig:MCillustration} shows MC results for a random circulant matrix \eqref{circulant matrix} with $p=100, n=1,000$ (top panel), and for $p=1,000, n=10,000$ (bottom panel). Results for a random  Toeplitz matrix \eqref{toeplitz matrix} are very similar, and we do not report them here. The left panel plots the average over $10,000$ MC replications of the eigenvalues of $\bS$, starting from the largest (order number 1 on the horizontal axis) and ending with the smallest (order number $p$). The right panel shows the histogram of the $10,000$ MC realizations of the smallest eigenvalue. 

We see that the smallest eigenvalue, $\lambda_p$, remains relatively insensitive to the ten-fold increase in the values of $p$ and $n$. The histogram of its MC replications concentrates and rather sluggishly shifts to the left. In fact, the left boundary of the histogram (more precisely, the smallest MC replication of $\lambda_p$) \textit{increases} from $0.1582$ to $0.1595$ when the values of $p$ and $n$ increase in our simulation experiment. Arguably, one may easily misinterpret such MC results as suggesting a convergence of the smallest eigenvalue to some positive value. 

Interestingly, the Monte Carlo analysis in \cite{bose10} is employed to illustrate their theoretical findings on the existence of the limiting spectral distribution of $\bS$. (See Theorems 1 and 4 in \cite{bose10} for the Toeplitz and circulant cases, respectively.) Their results suggest that, for the Toeplitz case with $c = 1/3$, the support of the limiting distribution is bounded away from zero. Inspired by this evidence, \cite{bose10} pose a question regarding the positivity of the infimum of the support and its precise value.

Our finding that $\lambda_p$ converges to zero casts doubt on the hypothesis that the support of the limiting spectral distribution is separated from zero when $c<1$. However, it is not inconsistent with this hypothesis. Conceivably, $o(p)$ of the relatively small eigenvalues could remain below a positive lower bound of the support of the limiting spectral distribution, so that the smallest eigenvalue does converge to zero and the hypothesis still holds. We leave proving or rejecting the hypothesis for future study.

There have been many previous studies of the asymptotic behavior of the singular values of $\bX/\sqrt{n}$ (square roots of the eigenvalues of $\bS$) for random Toeplitz and circulant\footnote{For a book-long treatment of ``patterned'' random matrices, including Toeplitz and circulant matrices, see \cite{bose18}.} $\bX$ under various assumptions on the generating sequence $x_i$. These studies can be roughly categorized into four strands. 

The first, most numerous, strand consists of the research on the asymptotic behavior of the empirical distribution of the singular values. For some, but not all, important works following this direction see \cite{bose02}, \cite{Bryc06}, \cite{hammonf05},  \cite{bose10},  \cite{bose23}, and \cite{shamb23}.
The second strand focuses on the fluctuations of the linear spectral statistics, that is, of the averages of some regular test functions of the singular values. We would like to mention, not pretending to be thorough, \cite{chatterjee09}, \cite{liu12}, and \cite{Adhikari17, Adhikari18}. The third group is represented by studies of the largest singular value. This group includes, but is not limited to \cite{meckes07}, \cite{bose07}, \cite{adamczak10,adamczak13}, and \cite{sen13}.

The most relevant for this paper is the fourth group, which constitutes the previous literature on the smallest eigenvalue of $\bS$ for the random Toeplitz and circulant $\bX$. This literature has been focusing on the square case $p=n$, hence $c=1$. \cite{pan15} derives probabilistic bounds on the smallest eigenvalue of $\bS$ for i.i.d.~Gaussian $x_i$. When $x_i$ are i.i.d.~and satisfy the Lyapunov condition, but not necessarily Gaussian, \cite{barrera22} proves that the smallest singular value of the circulant $\bX$ (which equals the square root of the product of $n$ and the smallest eigenvalue of $\bS$ for the circulant case) is asymptotically distributed according to the Rayleigh distribution with c.d.f. $R(x)=1-\exp\{-x^2/2\}, x\geq 0$. This, of course, implies that the smallest eigenvalue of $\bS$ itself converges to zero in probability. 

Note that when $p=n$, the smallest eigenvalue of $\bS$ converges to zero even in the standard case of sample covariance matrices when all elements of $\bX$ are i.i.d. For readers interested in the second-order asymptotic behavior of the smallest eigenvalue in this scenario, we refer to \cite{tao10}. Other notable works on the smallest eigenvalue of $\bS$ for square random $\bX$ under various assumptions on its entries are summarized in the introduction and Remark 1.5 of \cite{barrera22}.

To our knowledge, this is the first paper to analyze the asymptotic behavior of $\lambda_p$ for rectangular Toeplitz or circulant $\bX$ with $c<1$.

The remainder of the paper is organized as follows. The next section presents our main results (Theorems \ref{thm: exponential}-\ref{thm: small T/p}) and outlines their proofs. Section \ref{sec: MC} reports results of the Monte Carlo analysis.  Section \ref{sec: proof details} contains further details of the proofs of Theorems \ref{thm: exponential}-\ref{thm: small T/p}. Section \ref{sec: conclusion} concludes. Proofs of technical lemmas are given in the Appendix.

\section{Main results and proof outlines}\label{sec: main results}
In this paper, we consider the asymptotic regime where $p,n\rightarrow\infty$. We assume that $n:=n(p)$ is a function of $p$ such that as $p\rightarrow\infty$, $p/n\rightarrow c\in(0,1]$. We will abbreviate this as $p,n\rightarrow_c\infty$.
Recall that $\lambda_p$ denotes the smallest eigenvalue of $\bS=\frac{1}{n}\bX\bX^\top$, where $\bX$ is either a random Toeplitz matrix \eqref{toeplitz matrix} or a random ``rectangular circulant'' matrix \eqref{circulant matrix}. In what follows, we will omit the quotation marks around the term ``rectangular circulant'', which will always correspond to matrices of type \eqref{circulant matrix}.

\begin{theorem}
\label{thm: exponential}
    Suppose that $p,n\rightarrow_c\infty$, where $c\in(0,1]$. Then, for both Toeplitz and circulant cases, there exists a constant $\beta>0$ that may depend on $c$, such that
    \begin{equation}
    \label{eq: thm1} \lambda_p=o_{\mathrm{P}}\left(\exp\left\{-\beta\log^{1/3}(p)\right\}\right).
    \end{equation}
\end{theorem}

 Note that for any rectangular $p\times n$ matrix $\mathbf{X}$ with $p>n$, $\lambda_p=0$. Therefore, the convergence of $\lambda_p$ to zero is trivial for $c>1$. Hence, definitions \eqref{toeplitz matrix} and \eqref{circulant matrix} of $\bX$ as well as all our results focus on the non-trivial cases $c\in(0,1]$.
The following theorem shows that $\lambda_p$ converges to zero not only in probability, but also in expectation.

\begin{theorem}
\label{thm: expectation}Under assumptions of Theorem \ref{thm: exponential}, for every $\kappa\geq 1$, $\mathbb{E}(\lambda_p)^\kappa\rightarrow 0$.
\end{theorem}

Function $\exp\left\{-\beta\log^{1/3}(p)\right\}$ in \eqref{eq: thm1} converges to zero faster than $\log^{-\alpha}p$ but slower than $p^{-\alpha}$ for any $\alpha>0$. Hence, Theorem \ref{thm: exponential} implies that $\lambda_p$ converges to zero in probability faster than any inverse power of logarithm. Of course, the theorem does not give us any upper bound on the convergence rate because even the function which is identically equal to zero (instant convergence) can be classified as $o_{\mathrm{P}}\left(\exp\left\{-\beta\log^{1/3}(p)\right\}\right)$.

We were able to establish some polynomial upper bounds on the rate of convergence in probability only for the circulant case. The following theorems describe these upper bounds.
\begin{theorem}
\label{thm: lower bound}
    Suppose that $\bX$ is a random rectangular circulant matrix \eqref{circulant matrix}. Let $m\geq 2$ be a fixed integer, and   suppose that $p,n\rightarrow_c\infty$ with $c<\frac{1}{\pi(m-1)}$, or equivalently, that $\lim n/p>\pi(m-1)$. Then for any small $\epsilon>0$, we have
    \[
    \Pr\left\{\lambda_p> p^{-\frac{1}{m}-\epsilon}\right\}\rightarrow 1.
    \]

\end{theorem}

\begin{theorem}
\label{thm: small T/p}
    Under assumptions of Theorem \ref{thm: lower bound}, if $c\leq \frac{1-\delta}{m}$, or equivalently, $\lim n/p\geq m/(1-\delta)$, where $\delta>0$ is an  arbitrarily small fixed real number, then for any small $\epsilon>0$, we have
    \[
    \Pr\left\{\lambda_p>p^{-\frac{1}{2}-\frac{1}{2m}-\epsilon}\right\}\rightarrow 1.
    \]
    For $c\in [1/2,1]$, $\Pr\left\{\lambda_p>p^{-1-\epsilon}\right\}\rightarrow 1$.
\end{theorem}

Theorem \ref{thm: lower bound} implies that, for sufficiently small $c$, the rate of convergence of $\lambda_p$ to zero is arbitrarily polynomially slow. Specifically, $\lambda_p\overset{P}{\rightarrow}0$ slower than $p^{-\alpha}$, where $\alpha=\frac{\pi c}{\pi c+1}+\epsilon$ and $\epsilon>0$ is an arbitrarily small constant. However, the theorem is not providing any upper bound on the convergence rate for $c\geq \pi^{-1}$. 

Theorem \ref{thm: small T/p} uses different proof ideas to cover the case $\pi^{-1}\leq c<1/2$ at the expense of a much worse upper bound for smaller $c$. The statement of the theorem for $c\in [1/2,1]$ follows from the result of \cite{barrera22} about the convergence of $\sqrt{n\lambda_p}$ to the Rayleigh distribution in the square circulant case, mentioned in the Introduction.\vspace{2mm} 

\textbf{Proof outlines.} Here we introduce a setup and outline the proofs of the above theorems. Further details of the proofs can be found in Sections \ref{sec: thm1}-\ref{sec: thm4}. 

Consider a $p \times n$ rectangular circulant matrix $\bX$ as defined in \eqref{circulant matrix}. Let $\widetilde{\bX}$  be a Toeplitz matrix consisting of the last $\widetilde{n}:= n-p+1$ columns of $\bX$. It follows that $\bX \bX^\top - \widetilde \bX \widetilde \bX^\top \geq 0$, implying that the smallest eigenvalue of $\widetilde{\bX}\widetilde{\bX}^{\top}$ is less than or equal to the smallest eigenvalue of $\bX\bX^{\top}$. Therefore, to establish Theorems \ref{thm: exponential} and \ref{thm: expectation} for $p\times\tilde{n}$ Toeplitz matrices, it suffices to prove them for $p\times n$ circulant matrices. Consequently, throughout our proofs, we assume that $\bX$ is a rectangular circulant matrix.\vspace{1mm}

\textbf{About Theorem \ref{thm: exponential}:} The key element of the proof of Theorem \ref{thm: exponential} is an upper bound on $\lambda_p$, expressed in the form of a weighted average of the periodogram of $x_0,\dots,x_{n-1}$. The variational characterization of the smallest eigenvalue yields
\begin{equation}
\label{variational}
\lambda_p\leq \frac{1}{n}\bu^\ast \bX\bX^\top\bu
\end{equation}
for any complex $p$-dimensional vector $\bu:=(u_0,u_1,\dots,u_{p-1})^\top$ with Euclidean norm 1. Let $\bX_C$ be the $n\times n$ circulant matrix with the first $p$ rows constituting $\bX$, and let $\bu_C$ be an $n$-dimensional vector $\bu_C:=(u_0,u_1,\dots,u_{p-1},0,\dots,0)^\top$. Then obviously \begin{equation}
    \label{dilation} \bu^\ast\bX\bX^\top\bu=\bu_C^\ast\bX_C\bX^\top_C\bu_C.
\end{equation}

Since $\bX_C^\top$ is an $n\times n$ real circulant matrix, it admits the spectral decomposition\footnote{See e.g. \cite{golub96}, p. 202.} 
\begin{equation}
\label{spectral circulant}
\bX_C^\top=\sqrt{n}\cF^\ast\operatorname{diag}\left(\cF \bx\right) \cF,
\end{equation}
where $\cF$ is the (scaled) unitary $n\times n$ Discrete Fourier Transform (DFT) matrix with elements\footnote{We start indexing form $0$, so that the element in the first row and the first column of $\cF$ is denoted as $\cF_{0,0}$}  $\cF_{s,t}=\exp\{\frac{2\pi\rmi}{n}st\}/\sqrt{n}$, and 
$\mathbf{x}:=(x_0,x_{1},x_{2},\dots,x_{n-1})^\top$ is the first column of $\bX_C^\top$. 
Let $\by:= \cF\bx$ be the (scaled) DFT of $\bx$ and let $\mb{P}_\bu := \sqrt{n}\cF\bu_C$ be the DFT of $\bu_C$. Then using \eqref{spectral circulant} and \eqref{dilation} in \eqref{variational}, we obtain
\begin{equation}
\label{periodogam bound}
\lambda_p\leq \frac{1}{n}\bu_C^\ast\bX_C\bX^\top_C\bu_C=\frac{1}{n}\|\by\odot \mb{P}_\bu\|^2=\frac{1}{n}\sum_{s=0}^{n-1}|y_s|^2|P_{\bu s}|^2,
\end{equation}
where $\odot$ denotes the element-wise Hadamard product, and $y_s,P_{\bu s}$ are the $s$-th elements of $\mb{y}$ and $\mb{P}_\bu$.
In the proof, we will seek a $\mb{P}_\bu$ that ensures the right-hand side of \eqref{periodogam bound} is small with high probability.

In the sum in \eqref{periodogam bound}, $|y_s|^2$ can be interpreted as the periodogram of $\bx$ at the fundamental frequencies $\omega_s:= 2\pi s/n$, while  $P_{\bu s}=P_{\bu}\left(e^{\mathrm{i}\omega_s}\right)$, where
\[
P_{\bu}(z):=\sum_{j=0}^{p-1}u_j z^j
\]
is a polynomial with the coefficients given by the elements of $\bu$. Hence, \eqref{periodogam bound} bounds $\lambda_p$ by a weighted average of the periodogram of $\bx$ with weights determined by the magnitude of the polynomial $P_{\bu}(z)$ of degree $p-1$ on a grid of size $n$ on the unit circle.  

Since $\mathcal{F}$ is unitary and $\bx$ is normally distributed, $\by= \mathcal{F}\bx$ has a complex normal distribution.  Precisely,
\begin{equation}
\label{y definition}
\mathbf{y}=\begin{cases}
\begin{pmatrix}y_0&y_1&\dots&y_{(n-1)/2}&\overline{y_{(n-1)/2}}&\dots&\overline{y_1}\end{pmatrix}^\top&\text{ if }n\text{ is odd}\\
\begin{pmatrix}
y_0&y_1&\dots&y_{n/2-1}&y_{n/2}&\overline{y_{n/2-1}}&\dots&\overline{y_1}
\end{pmatrix}^\top&\text{ if }n\text{ is even}
\end{cases},
\end{equation}
where $y_0,y_{n/2}$ are i.i.d.~$\mathcal{N}(0,1)$ and $\{y_i, 0<i<n/2\}$ are independent from them i.i.d.~$\mathcal{N}_{\mathbb{C}}(0,1)$ (that is, random variables of form $a_1+\mathrm{i}a_2$, where $a_1$ and $a_2$ are independent $\mathcal{N}(0,1/2)$). In particular, for any $0<s<t<n/2$, $|y_s|^2$ and $|y_t|^2$ are independent random variables having exponential distribution $\mathrm{Exp}(1)$. 
   We also note that the facts that $\mathcal{F}$ is unitary and $\|\bu\|=1$ imply that $\frac{1}{n}\sum_{s=0}^{n-1}|P_{\bu s}|^2=1$.
   
   If we choose $\bu$ so that $|P_{\bu s}|^2$ for all $s$ are roughly of the same order of magnitude, the law of large numbers would imply that the right hand side of \eqref{periodogam bound} is $O_P(1)$, which would not be helpful as a bound for $\lambda_p$. Hence, we will aim at choosing $\bu$ so that \textcolor{black}{ $|P_{\bu s}|^2$ very quickly decays when $\min\{s,n-s\}$ grows.}

   \textcolor{black}{As follows from the DFT interpretation of $\mb{P}_\bu$, if instead of vector $\bu$ with coordinates $u_s$, $s=0,\dots,p-1$ we consider vector $\bu_\tau$ with coordinates $\exp\{-\omega_\tau s \rmi\}u_s$, $s=0,\dots,p-1$, the corresponding $P_{\bu_\tau, s}$ would satisfy the identity
   $P_{\bu_\tau, s}=P_{\bu,(s-\tau)\operatorname{mod} n}$.
   In what follows, we will omit notation $\operatorname{mod}n$ from any coordinate indices of $n$-dimensional vectors. For example, we simply write 
\[
P_{\bu_\tau, s}=P_{\bu, s-\tau}.
\]
Pickings several integers $\tau_1,\dots,\tau_K$, defining $\bu_k:=\bu_{\tau_k}$, and using the above identity, we obtain the following stronger version of bound \eqref{periodogam bound}: 
\begin{equation}
\label{stronger periodogram bound}
\lambda_p\leq \min_{k=1,\dots,K}\frac{1}{n}\sum_{s=0}^{n-1}|y_s|^2|P_{\bu,s-\tau_k}|^2.
\end{equation}
Since by construction the coordinates of $\mb{P}_\bu$ will be quickly decaying, picking integers $\tau_1,\dots,\tau_K$ so that they are sufficiently far apart will ensure that the different sums under the minimum are almost independent, hence it will be relatively easy to evaluate the minimum. }
   

Recall that $\bP_\bu$ is the DFT of vector $\bu_C$, that ``adds''  $n-p$ zero coordinates to $\bu$. If we interpret this vector as a signal ``localized in time'' in the sense that only the first $p$ observations (corresponding to the entries of vector $\bu$) are non-zero, then \textcolor{black}{the task of choosing $\bu$ so that the coordinates of $\mb{P}_\bu$ quickly decay is similar} to the problem of choosing a time-localized signal that has the DFT which is well concentrated in the frequency domain. There exists a large literature on this topic, see \cite{slepian78}, \cite{donoho89}, \cite{barnett22} and references therein. 

In Section \ref{sec: thm1}, we follow \cite{barnett22} and consider the Gaussian Fourier pair,\footnote{\cite{barnett22} also considers the so-called Kaiser-Bessel Fourier pair. It turns out that in our setting, it leads to the same results, so we do not consider it here.} which is well localized in both time and frequency domain. This choice will imply tight bounds on $|P_{\bu, s-\tau_k}|^2$ in \eqref{stronger periodogram bound} yielding the convergence of $\lambda_p$ to zero with the speed stated in Theorem \ref{thm: exponential}.

\textbf{About Theorem \ref{thm: expectation}:} Theorem \ref{thm: exponential} and the continuous mapping theorem imply that $\left(\lambda_p\right)^\kappa\overset{\mathrm{P}}\rightarrow 0$ as $p \to \infty$.
To prove Theorem \ref{thm: expectation}, it is sufficient to show that $\E (\lambda_p)^{\kappa + \eps} \leq C < \infty$ for some $\eps > 0$, where $C$ is a constant that does not depend on $p$. Indeed, this implies that $(\lambda_p)^\kappa$ is  uniformly integrable and therefore the convergence in probability implies the convergence in expectation (Theorem 5.5.4 in \cite{gut2013}). In Section \ref{sec: thm2}, we prove that 
\begin{equation}
\label{to prove for thm2}
  \E \lambda_p^\kappa \leq C_{\kappa, c} < \infty\qquad\text{for every}\quad \kappa \geq 1  
\end{equation} 
with a constant $C_{\kappa,c}$ that may depend on $\kappa$ and $c:=\lim p/n$. \textcolor{black}{Our detailed proof in Section \ref{sec: thm2} relies on the inequality \eqref{periodogam bound} and inequalities for sub-exponential random variables (see sections 2.7 and 2.8 of \cite{vershynin2018}).}\vspace{1mm}

\textbf{About Theorem \ref{thm: lower bound}:} The starting point of our proof of Theorem \ref{thm: lower bound} is the fact that 
\[
\lambda_p=\sum_{j=0}^{n-1}|y_j|^2\left(\frac{1}{n}|P_\bu(e^{\rmi\omega_j})|^2\right)
\]
\textit{for some} real $\bu$ with $\|\bu\|=1$. We aim to show that the right-hand side of the above equality is greater than $p^{-\frac{1}{m}-\epsilon}$ \textit{for all} real $\bu$ with $\|\bu\|=1$, for all sufficiently large $p$ on a sequence of events $A_p$ such that $\Pr A_p\rightarrow 1$.

\textcolor{black}{Note that $\frac{1}{n}|P_\bu(e^{\rmi\omega})|^2$ is a cosine trigonometric polynomial of order $p-1$. We denote it as $T_{p-1}(\omega)$, suppressing its dependence on $\bu$. Since $\|\bu\|=1$, we have
\begin{equation}
\label{trigonometric normalization}
\int_0^{2\pi}|T_{p-1}(\omega)|\mathrm{d}\omega =\int_0^{2\pi}\frac{1}{n}|P_\bu(e^{\rmi \omega})|^2\mathrm{d}\omega=\frac{2\pi}{n}
\end{equation}
and, assuming for concreteness that $n$ is odd (the case of even $n$ requires a similar analysis, and we omit it) and letting $N:=(n+1)/2$,
\begin{equation}
    \label{discrete trigonometric normalization}
    T_{p-1}(0)+2\sum_{j=1}^{N-1}T_{p-1}(\omega_j)=\sum_{j=0}^{n-1}\frac{1}{n}|P_\bu(e^{\rmi\omega_j})|^2=1.
\end{equation}
Further, define
\[
\xi_0:=|y_0|^2,\qquad \xi_j:=|y_j|^2+|y_{n-j}|^2=2|y_j|^2, \quad j=1,\dots,N-1.
\]
With this notation, we have 
\begin{eqnarray} \notag
\lambda_p= \sum_{j=0}^{N-1}\xi_j T_{p-1}(\omega_j)
&=&\xi_0T_{p-1}(0)+\sum_{j=1}^{N-1}\xi_{(j)} T_{p-1}(\omega_{\sigma(j)})\\\label{lambdap for thm3}
&\geq& \frac{\xi_0}{2}T_{p-1}(0)+\xi_{(r)}\sum_{j=r+1}^{N-1} T_{p-1}(\omega_{\sigma(j)}),
\end{eqnarray}
where $\xi_{(j)}$ with $j\geq 1$ is the $j$-th order statistic for the sequence $\xi_1,\dots,\xi_{N-1}$,  $\sigma$ is the (random) permutation of the indices such that $\xi_{(j)}=\xi_{\sigma(j)}$, and $r:=\lfloor p^{1-1/m-\epsilon/2}\rfloor$. }

\textcolor{black}{We will assume that $\min\{\xi_0,\xi_{(r)}\}>r/N$ on events $A_p$. We can make such an assumption because  
\begin{equation}
\label{to show in sec 5}
\Pr\left(\min\{\xi_0,\xi_{(r)}\}>r/N\right)\rightarrow 1
\end{equation}
as $p\rightarrow\infty$. Section \ref{sec: thm3} contains a quick proof of \eqref{to show in sec 5}. Therefore, \eqref{lambdap for thm3} implies that on $A_p$,
\begin{equation}
\label{referrred to next}
\lambda_p>\frac{r}{N} \left(\frac{1}{2}T_{p-1}(0)+\sum_{j=r+1}^{N-1}T_{p-1}(\omega_{\sigma(j)})\right)=:\frac{r}{N}S_1.
\end{equation}
Since $\frac{r}{N}\gg p^{-1/m-\epsilon}$, Theorem \ref{thm: lower bound} would follow if we show that $S_1$ is bounded away from zero over all non-negative trigonometric polynomials $T_{p-1}$, satisfying \eqref{trigonometric normalization} and \eqref{discrete trigonometric normalization}.}

\textcolor{black}{The main idea of our proof is that a small value of $S_1$ yields a large total variation of $T_{p-1}(\omega)$ on $[0,2\pi]$, with high probability. Indeed, normalization \eqref{discrete trigonometric normalization} implies that, if values of $T_{p-1}(\omega_{\sigma(j)})$ are small for $j>r$ (so that $S_1$ is small), they must be large for some $j\leq r$. On the other hand, positions of $\omega_{\sigma(j)}$ with $j\leq r$ and with $j>r$ in the interval $[0,\pi]$ are randomly mixed. This leads to a large total variation of $T_{p-1}(\omega)$ (see Section \ref{sec: thm3} for details).}

\textcolor{black}{As we further show in Section \ref{sec: thm3}, such a large total variation of $T_{p-1}(\omega)$ contradicts the integral version of the Bernstein inequality for the derivative of trigonometric polynomials (see e.g.~Theorem 2.5 in Chapter 4 of \cite{devore93}). Hence, $S_1$ cannot be small, which concludes our proof of Theorem \ref{thm: lower bound}.}
\vspace{1mm}

\textbf{About Theorem \ref{thm: small T/p}: } \textcolor{black}{Unfortunately, the Bernstein inequality for trigonometric polynomials does not allow us to obtain a non-trivial upper bound on the convergence rate of $\lambda_p$ for $p, n$ such that $\lim n/p\leq\pi$. To cover this case we use Lemma \ref{lem: Stechkin 1} (due to \cite{stechkin48}, see Section \ref{sec: thm4}) showing that trigonometric polynomials cannot decay too fast near the point of their maximum.}

\textcolor{black}{Using an inequality similar to \eqref{referrred to next}, but with $r_1\approx \sqrt{r}$ replacing $r$, and assuming that the right-hand side of such an inequality is small,  Section \ref{sec: thm4} shows that the maximum of the corresponding trigonometric polynomial on $[0,\pi]$ must be relatively large. Since, by Lemma \ref{lem: Stechkin 1}, the  trigonometric polynomial must remain large in a relatively large neighborhood of the maximum point, such a neighborhood will include points $\omega_{\sigma{(j)}}$ with $j>r$, with high probability. The original inequality \eqref{referrred to next} then would imply that $\lambda_p$ cannot be too small.}

\textcolor{black}{Balancing the lower bounds on $\lambda_p$ provided by the modified and the original inequality \eqref{referrred to next} gives us the upper bound on the rate of the convergence of $\lambda_p$ to zero, stated by Theorem \ref{thm: small T/p}.}

\section{Monte Carlo}\label{sec: MC}

Theorems \ref{thm: exponential}-\ref{thm: small T/p} leave a considerable gap between the lower and upper bounds on the rate of convergence of $\lambda_p$ to zero. In this section, we perform Monte Carlo (MC) analysis to get some idea about the actual rate of the convergence. 
We simulate $p\times n$ rectangular circulant and Toeplitz matrices $\bX$, and compute the smallest eigenvalue $\lambda_p$ of $\bX\bX^\top/n$. We make 10,000 MC replications for four different ratios: $n/p=2,3,5,$ and $10$, and seven different values of $p$: $100,\dots,700$. Our simulations were done in MATLAB on a basic laptop computer. Therefore, allowing for $p>700$ would have been too time consuming. 

\begin{figure}[h]
\centering
\includegraphics[]{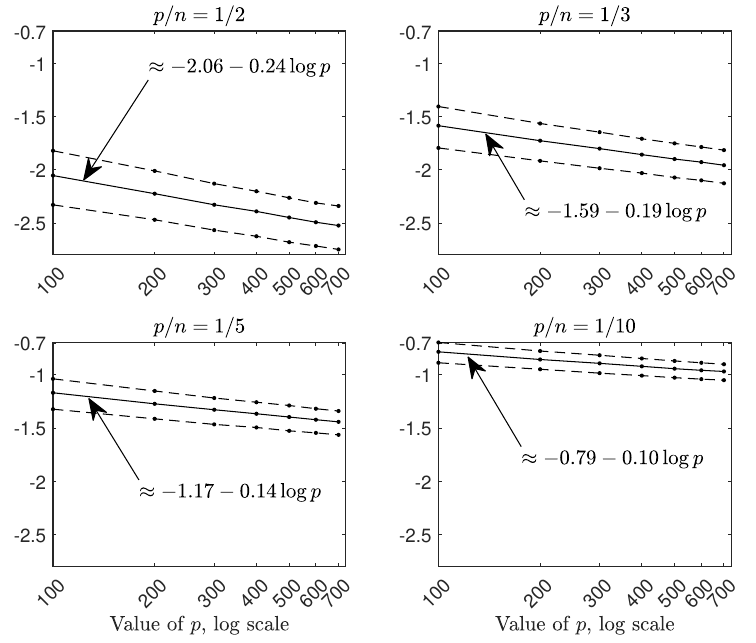}
\caption{25,50, and 75 percentiles of the Monte Carlo distribution of $\log \lambda_p$ for $p\times n$ rectangular circulant $\bX$. Based on 10,000 MC replications. Equations reported in the figure correspond to the ordinary least squares estimates of the median line based on seven observations (dot markers).}
\label{fig:circ speed log}
\end{figure}

Figure \ref{fig:circ speed log} reports the logarithm of $k=25\%,  50\%$, and $75\%$ quantiles, $Q_k(\lambda_p)$, of the Monte Carlo distribution of $\lambda_p$ for rectangular circulant matrix $\bX$. We plot $\log Q_k(\lambda_p)$ against $\log p$. The solid lines correspond to the median ($k=50\%$), the dashed line above it -- to $k=75\%$, and the dashed line below it -- to $k=25\%$. All the graphs look linear for all values of $p/n=1/2, 1/3, 1/5,$ and $1/10$. The lines representing the graphs for different percentiles $k$ are almost parallel. These lines become flatter as $p/n$ decreases. At the same time, the intersections of the lines with the vertical axis become closer to zero.
This suggests that the model
\begin{equation}
\label{quantile model}
Q_k(\lambda_p)=\alpha_{c,k} p^{-\beta_c},
\end{equation}
with some $\alpha_{c,k}>0$, increasing in both $k$ and $c=\lim p/n$, and $\beta_c>0$, decreasing in $c$,  might be approximating the asymptotic behavior of the quantiles of $\lambda_p$ reasonably well. If so, we would expect a polynomial rate of the convergence of $\lambda_p$ to zero. Furthermore, this rate becomes slower ($\beta_c$ decreases) as $c$ decreases. 

We estimate the slope $-\beta_c$ of the solid lines (log of the MC median of $\lambda_p$)  reported in Figure \ref{fig:circ speed log} by ordinary least squares (OLS)  regression based on seven observations, corresponding to $p=100,\dots, 700$. The estimates are reported next to the graphs. 
It is interesting to compare the estimates of $\beta_c$, $\hat{\beta}_c$, with the upper bounds on the convergence rates established in Theorem \ref{thm: lower bound}. 
For $p/n=1/5$, $\hat{\beta}_c=0.14$, whereas the smallest upper bound provided by Theorem \ref{thm: lower bound} for $p/n=1/5$ equals $(5/\pi+1)^{-1}+\epsilon\approx 0.39$, which is approximately $2.8$ times larger. For $p/n=1/10$, $\hat{\beta}_c=0.10$, whereas the smallest upper bound provided by Theorem \ref{thm: lower bound} for $p/n=10$  is $(10/\pi+1)^{-1}+\epsilon\approx 0.24$, which is $2.4$ times larger. Hence, the validity of Theorem \ref{thm: lower bound} is supported by our MC results. At the same time, these results suggest that the bounds described by Theorem \ref{thm: lower bound} are far from optimal.

For $p/n>1/\pi$, Theorem \ref{thm: lower bound} does not provide us with any bounds, so we need to use Theorem \ref{thm: small T/p} instead. For $p/n=1/2$, $\hat{\beta}_c=0.24$, whereas the smallest upper bound provided by Theorem \ref{thm: small T/p} for $p/n=1/2$ is approximately $3/4$, which is about $3.1$ times larger. For $p/n=1/3$, $\hat{\beta}_c=0.19$, whereas the smallest upper bound provided by Theorem \ref{thm: small T/p} for $p/n=1/3$ is approximately $2/3$, which is about $3.5$ times larger. Hence, again, MC supports the theoretical findings, but the bounds described by Theorem \ref{thm: small T/p} might be suboptimal.

\begin{figure}[h]
\centering
\includegraphics[]{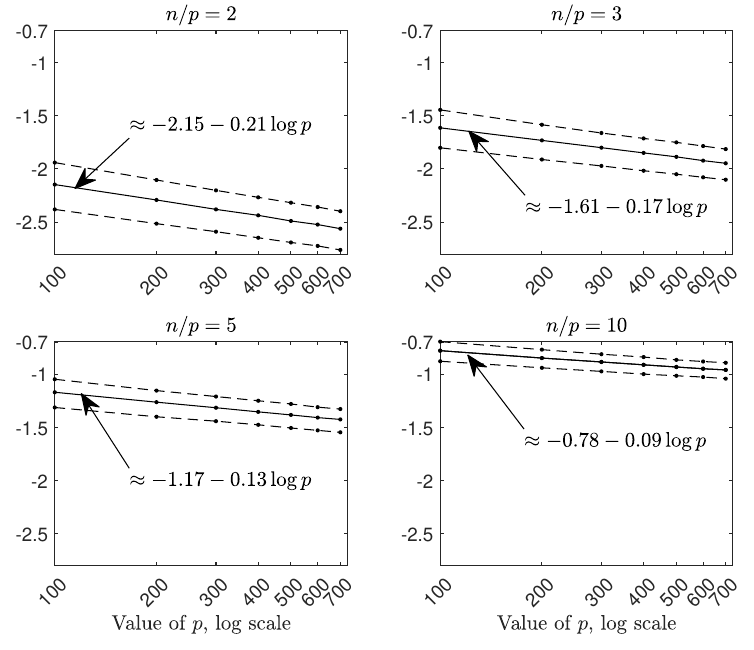}
\caption{25,50, and 75 percentiles of the Monte Carlo distribution of $\log \lambda_p$ for $p\times n$ rectangular Toeplitz $\bX$. Based on 10,000 MC replications. Equations reported in the figure correspond to the ordinary least squares estimates of the median line based on seven observations (dot markers).}
\label{fig:toep speed log}
\end{figure}

 Figure \ref{fig:toep speed log} is the equivalent of Figure \ref{fig:circ speed log} for Toeplitz $\bX$. The figures are very similar, and the Ordinary Least Squares (``OLS'') estimates of the slopes $-\beta_c$ differ only slightly (being smaller in absolute value) from those for rectangular circulant $\bX$. This suggests that extrapolating the results of Theorems \ref{thm: lower bound} and \ref{thm: small T/p} from circulant to Toeplitz matrices would not lead to a substantial error.

\section{Proof details}\label{sec: proof details}
This section contains details of the proofs of Theorems \ref{thm: exponential}-\ref{thm: small T/p}, outlined in Section \ref{sec: main results}.

\subsection{Proof details for Theorem \ref{thm: exponential} }\label{sec: thm1}
\textcolor{black}{As explained in Section \ref{sec: main results}, the proof is based on the inequality \eqref{stronger periodogram bound}. We start from making a choice of vector $\bu$.} We will need the following lemma. Its proof (see Appendix) is only slightly different from the proof of Proposition 6 in \cite{barnett22}.

\begin{lemma}[Proposition 6, \cite{barnett22}]
\label{lem: barnett}
    Let $\sigma^2_p$ be a positive constant that depends only on $p$, and let $\mathbf{w}=(w_0,...,w_{n-1})^\top$ be the so-called periodized Gaussian vector with
    \begin{equation}
    \label{barnett direct}
    w_j=(2\pi\sigma_p^2)^{-1/2}\sum_{s\in\mathbb{Z}}e^{-\frac{1}{2}(j-\lceil p/2\rceil +sn)^2/\sigma^2_p},\quad j=0,\dots,n-1,
    \end{equation}
    where $\lceil x\rceil$ denotes the smallest integer no smaller than $x$. 
    Then the DFT of $\mathbf{w}$, that is the vector $\hat{\mathbf{w}}=\sqrt{n}\mathcal{F}\mathbf{w}$, is also a periodized Gaussian vector with components
    \begin{equation}
    \label{branett transformed}
    \hat{w}_k=e^{\frac{2\pi\rmi}{n}\lceil p/2\rceil k}\sum_{m\in\mathbb{Z}}e^{-2\left(\frac{\pi\sigma_p}{n}\right)^2(k+mn)^2}, \quad k=0,\dots,n-1.
    \end{equation}
\end{lemma}


\textcolor{black}{We define $\bu$ as a truncated and normalized version of $\mb{w}$:}
\[
\bu:=\frac{1}{\sqrt{\sum_{s=0}^{p-1}w_s^2}}\left(w_0,\dots,w_{p-1}\right)^\top.
\]
The following lemma uses Lemma \ref{lem: barnett} to show that \textcolor{black}{the coordinates of the corresponding vector $\mb{P}_\bu$ quickly decay.} 
A proof of the lemma is in the Appendix.
\begin{lemma}
\label{lem: exponential d}
Suppose that $p,n\rightarrow_c\infty$ \textcolor{black}{with $c\in(0,1]$ and $p\leq n$}, and that $\sigma_p\rightarrow\infty$ so that $\sigma_p=o(p)$ and $\sqrt{p}=o(\sigma_p)$. Then for all sufficiently large $p$, we have
\begin{equation}
\label{lemma 2 ineq}
\frac{1}{n}|P_{\bu s}|^2\leq \frac{2^6\sigma_p}{n}e^{-4\left(\frac{\pi\sigma_p}{n }\right)^2s^2}+\frac{2^8 \sigma_p^3}{n p^2}e^{-\frac{p^2}{4\sigma_p^2}},\qquad s\in [0, n/2].
\end{equation} 
Furthermore, $|P_{\bu s}|=|P_{\bu,n-s}|=|P_{\bu,-s}|$.
\end{lemma}

\textcolor{black}{Changing the summation index in \eqref{stronger periodogram bound} and recalling that we understand indices modulo $n$, we obtain}
\begin{equation}
\label{strong periodogram bound 1}
\lambda_p\leq \min_{k=1,\dots,K}\frac{1}{n}\sum_{s=0}^{n-1}|y_s|^2|P_{\bu, s-\tau_k }|^2=\min_{k=1,\dots,K}\frac{1}{n}\sum_{s=0}^{n-1}|y_{s+\tau_k}|^2|P_{\bu s }|^2.
\end{equation}
\textcolor{black}{Let us choose the integers $\tau_k$, $k=1,\dots,K$ as follows. Consider } slowly increasing integers $b_p=o(p)$ and  let $K:=\lfloor \frac{n}{4b_p+2}\rfloor-1$, where $\lfloor\cdot\rfloor$ denotes the integer part of a real number. Let
\begin{equation}
\label{thetak definition}
\tau_k:=k(2b_p+1),\quad k=1,2,\dots,K.
\end{equation}




For any $k=1,\dots, K$, define
\begin{eqnarray*}
Y_k&:=&\frac{1}{n}\sum_{s=-b_p}^{b_p}|y_{s+\tau_k}|^2|P_{\bu s}|^2\quad\text{and}\\
R_k&:=&\frac{1}{n}\sum_{s=0}^{n-1}|y_{s+\tau_k}|^2|P_{\bu s }|^2-Y_k.
\end{eqnarray*}
Note that \eqref{y definition} and our definitions of $b_p$ and $K$ imply that $Y_k$, $k=1,\dots,K$, are independent and identically distributed. Furthermore, $R_k$ depend only on $P_{\bu s}$ with $b_p<s<n-b_p$, and hence, by Lemma \ref{lem: exponential d} they can be shown to be ``small'' (after carefully choosing $\sigma_p$ and $b_p$ sequences). 

Keeping this in mind, we use \eqref{strong periodogram bound 1} to obtain the following bound on the probability that \textcolor{black}{$\lambda_p\leq e^{-2\beta\log^{1/3} p}$ with  $\beta>0$ possibly depending on $c$. Let $\epsilon_p=\frac{1}{2}e^{-2\beta\log^{1/3} p}$, then we have} 
\begin{eqnarray}
\notag
\Pr\left( \lambda_p\leq e^{-2\beta\log^{1/3} p}\right)&\geq& \Pr\left(\min_{k=1,\dots,K}(Y_k+R_k)\leq 2\epsilon_p\right)\\ \notag
&\geq& \Pr\left(\min_{k=1,\dots,K}Y_k\leq \epsilon_p\text{ and } \max_{k=1,\dots,K}R_k\leq \epsilon_p\right)\\ \notag
&\geq&\Pr\left(\min_{k=1,\dots,K}Y_k\leq \epsilon_p\right)-\Pr\left( \max_{k=1,\dots,K}R_k> \epsilon_p\right)\\ \notag
&=&1-\left(\Pr\left(Y_k> \epsilon_p\right)\right)^K-\Pr\left( \max_{k=1,\dots,K}R_k> \epsilon_p\right)\\ \label{the main inequality}
&\geq&1-\left(1-\Pr\left(Y_k\leq \epsilon_p\right)\right)^K-\sum_{k=1}^K\Pr\left( R_k> \epsilon_p\right).
\end{eqnarray}
\textcolor{black}{The remaining two steps of the proof will establish}
 a lower bound on $\Pr\left(Y_k\leq \epsilon_p\right)$ and an upper bound on $\Pr\left( R_k> \epsilon_p\right)$ \textcolor{black}{that would imply that the right hand side of \eqref{the main inequality} converges to unity. Such a convergence yields
 \[
 \Pr\left(\frac{\lambda_p}{\exp\{-\beta\log^{1/3} p\}}\leq \exp\{-\beta\log^{1/3} p\}\right)\rightarrow 1,
 \]
 which implies that $\lambda_p=o_{\mathrm{P}}\left(\exp\{-\beta\log^{1/3} p\}\right)$ as stated by Theorem \ref{thm: exponential}.}\vspace{2mm}

\textbf{Step 1: establishing a lower bound on $\Pr\left(Y_k\leq \epsilon_p\right)$. }
Using \eqref{y definition} in the definition of $Y_k$ yields\footnote{\textcolor{black}{We remind the reader that $\mathrm{Exp}(1)$ can be interpreted as $1/2$ times chi-squared with two degrees of freedom.}}
\[
Y_k=\frac{1}{n}\sum_{s=-b_p}^{b_p}\frac{1}{2}\chi^2_s(2)|P_{\bu s}|^2,
\]
where $\chi^2_s(a)$, $s=-b_p,\dots, b_p$ are i.i.d.~random variables having chi-squared distribution with $a$ degrees of freedom. Denoting $\frac{1}{2n}|P_{\bu s}|^2$ as $d_s$ and recalling that by Lemma \ref{lem: exponential d} $|P_{\bu s}|=|P_{\bu,-s}|$, we represent $Y_k$ in the form
\[
Y_k=d_0\chi^2_0(2)+\sum_{s=1}^{b_p}d_s\chi^2_s(4)\leq \sum_{s=0}^{b_p}d_s\chi^2_s(4),
\]
where $\chi^2_s(4)$ is obtained from $\chi^2_s(2)$ by adding an independent chi-squared random variable with two degrees of freedom.

By Lemma \ref{lem: exponential d}, for all $0\leq s\leq n/2$ and all sufficiently large $p$, we have
\begin{equation}
\label{bound on d}
d_s\leq r_s+\Delta,\quad\text{where }r_s:=\frac{2^5\sigma_p}{n}e^{-4\left(\frac{\pi\sigma_p}{n}\right)^2s^2},\quad \Delta:=\frac{2^7\sigma_p^3}{np^2}e^{-\frac{p^2}{4\sigma_p^2}}.
\end{equation}
Therefore, we have
\[
Y_k\leq \sum_{s=0}^{\infty}r_s\chi^2_s(4)+\Delta \chi^2(4(b_p+1))=: Z_1+Z_2,
\]
where $\chi^2(4(b_p+1))$ is a random variable having chi-squared distribution with $4(b_p+1)$ degrees of freedom. Clearly,
\begin{equation}
\label{split Y into Zs}
\Pr(Y_k\leq\epsilon_p)\geq \Pr(Z_1\leq\epsilon_p/2)-\Pr(Z_2>\epsilon_p/2).
\end{equation}
We now establish a lower bound on $\Pr(Z_1\leq\epsilon_p/2)$ and an upper bound on $\Pr(Z_2>\epsilon_p/2)$.


First, we find a lower bound on $\Pr(Z_1\leq\epsilon_p/2)$. To simplify expressions to appear, we will introduce the following new notation:
\[
\alpha_p:= 2\pi\sigma_p/n,\qquad \tilde{\epsilon}_p:=\pi\epsilon_p/2^6,\qquad \tilde{Z}_1:=\sum_{s=0}^\infty\frac{1}{2}\alpha_pe^{-(\alpha_ps)^2}\chi_s^2(4),
\]
so that
\[
\Pr(Z_1\leq\epsilon_p/2)=\Pr(\tilde{Z}_1\leq \tilde{\epsilon}_p).
\]
To derive a lower bound on $\Pr(\tilde{Z}_1\leq\tilde{\epsilon}_p)$, we are going to use the upper and lower bounds on the probability of small Gaussian ellipsoids derived in section 2 of \cite{mayerwolf93}.

\textcolor{black}{That paper considers the random variable $z=\sum_{i=1}^\infty x_i^2/a_i^2$, where $a_i$, $i=1,2,\dots$ is a deterministic sequence with $\sum a_i^{-2}<\infty$ and $x_i$ are i.i.d.~standard normal random variables. It derives an upper bound $UB$ and a lower bound $LB$ on $\Pr\left(z<\epsilon\right)$ and establishes a simple inequality (their inequality (15)) that bounds the ratio $UB/LB$ from above. This simple inequality can be used to obtain a simple lower bound on $LB$ in terms of $UB$. We are going to use such a simple lower bound below.}

\textcolor{black}{Note that $\tilde{Z}_1$ can be represented in the form $\sum_{i=1}^\infty x_i^2/a_i^2$ with \[
a_{4s+1}^2=\dots=a_{4s+4}^2=2e^{(\alpha_p s)^2}/\alpha_p,\qquad s=0,1,\dots.
\]
\cite{mayerwolf93}'s upper bound on $\Pr(\tilde{Z}_1\leq\tilde{\epsilon}_p)$ is derived as follows.}
Since for any $t<1/2$,
\[
\mathbb{E}\exp\{t\chi_j^2(4)\}=(1-2t)^{-2},
\]
we have for any $\ell_p>0$,
\begin{equation}
\label{wolf-zeit ub}
\Pr(\tilde{Z}_1\leq\tilde{\epsilon}_p)\leq\mathbb{E}\left(\exp\left\{-\ell_p(\tilde{Z}_1-\tilde{\epsilon}_p)\right\}\right)=\exp\{\ell_p\tilde{\epsilon}_p\}\prod_{s=0}^\infty\left(1+\ell _p\alpha_p e^{-(\alpha_ps)^2}\right)^{-2}=:UB.
\end{equation}
\textcolor{black}{\cite{mayerwolf93} denote such a $UB$ as $UB(s_\epsilon)$, where their $s_\epsilon$ corresponds to our $\ell_p\tilde{\epsilon}_p$ and their $\epsilon$ corresponds to our $\tilde{\epsilon}_p$.} 

Now define $\delta(\ell_p,\tilde{\epsilon}_p)$ as
\begin{equation}
\label{delta ell epsilon}
\delta(\ell_p,\tilde{\epsilon}_p)=1-\frac{2}{\ell_p\tilde{\epsilon}_p}\sum_{s=0}^\infty\frac{1}{(\ell_p \alpha_p)^{-1}e^{(\alpha_ps)^2}+1}.
\end{equation}
\textcolor{black}{This is an equivalent of equation (11) in \cite{mayerwolf93}.} Suppose that $\ell_p$ is chosen so that $\delta(\ell_p,\tilde{\epsilon}_p)\in(0,1/2)$ (we verify that such a choice is possible below). Then, as follows from (15) of \cite{mayerwolf93}, \textcolor{black}{a lower bound $LB$ on $\Pr(\tilde{Z}_1\leq \tilde{\epsilon}_p)$ satisfies
\[
LB\geq UB \exp\{-2\ell_p\tilde{\epsilon}_p\delta(\ell_p,\tilde{\epsilon}_p)\}\left(1-\frac{1}{\ell_p\tilde{\epsilon}_p\delta^2(\ell_p,\tilde{\epsilon}_p)}\right).
\] The explicit form of $LB$ is of no interest to us here. Taking logarithm of both sides and using \eqref{wolf-zeit ub}, we obtain}
\begin{eqnarray}
\label{MWZ lower bound}
\log \Pr(\tilde{Z}_1\leq\tilde{\epsilon}_p)\geq \left(1-2\delta(\ell_p,\tilde{\epsilon}_p)\right)\ell_p\tilde{\epsilon}_p&-&2\sum_{s=0}^{\infty}\log \left(1+\ell_p \alpha_p e^{-(\alpha_ps)^2}\right)\\\notag
&+&\log\left(1-\frac{1}{\ell_p\tilde{\epsilon}_p\delta^2(\ell_p,\tilde{\epsilon}_p)}\right).
\end{eqnarray}
We will need the following lemma. Its proof is in the Appendix.
\begin{lemma}
    \label{lem: prelim}
    For any sequences $\ell_p>0$ and $\alpha_p>0$ such that $\ell_p\alpha_p\rightarrow\infty$ while $\alpha_p\rightarrow 0$ and $\alpha_p\log^{1/2}(\ell_p\alpha_p)\rightarrow 0$, we have
    \begin{eqnarray}
        \label{log asymptotics}
        \sum_{s=0}^\infty\log\left(1+\ell_p\alpha_pe^{-(\alpha_ps)^2}\right)&=&\frac{2}{3}\frac{\log^{3/2}(\ell_p\alpha_p)}{\alpha_p}(1+o(1)),\\\label{log der asymptotics}
        \sum_{s=0}^\infty\frac{1}{(\ell_p\alpha_p)^{-1}e^{(\alpha_ps)^2}+1}&=&\frac{\log^{1/2}(\ell_p\alpha_p)}{\alpha_p}(1+o(1)).
    \end{eqnarray}
    
\end{lemma}
Now 
define $\sigma_p$ and $\ell_p$ so that
\begin{equation}
\label{all def}
\alpha_p=\frac{2\pi\sigma_p}{n}=\left(\frac{2^5{\beta}^{3}}{\log p}\right)^{1/2},\quad\text{and}\quad \ell_p \alpha_p=\frac{3\log^{1/2} (1/\epsilon_p)}{\tilde{\epsilon}_p},
\end{equation}
\textcolor{black}{and recall that
\[
\tilde{\epsilon}_p=2^{-6}\pi\epsilon_p=2^{-7}\pi e^{-2\beta\log^{1/3} p}.
\]}
Note that the conditions of Lemma \ref{lem: prelim} are satisfied, so that \textcolor{black}{from \eqref{delta ell epsilon} and  \eqref{log der asymptotics},}
\begin{equation}
\label{delta less than 1/2}
\delta(\ell_p,\tilde{\epsilon}_p)=1-\frac{2\log^{1/2}(\ell_p\alpha_p)}{\ell_p\alpha_p\tilde{\epsilon}_p}(1+o(1))=\frac{1}{3}(1+o(1)),
\end{equation}
which belongs to $(0,1/2)$ for sufficiently large $p,n$, as required.

Further, the first two terms on the right hand side of \eqref{MWZ lower bound} satisfy
\begin{eqnarray*}
    (1-2\delta(\ell_p,\tilde{\epsilon}_p))\ell_p\tilde{\epsilon}_p&=&\frac{\log^{1/2} (1/\epsilon_p)}{\alpha_p}(1+o(1)),\\
    -2\sum_{s=0}^\infty\log\left(1+\ell_p \alpha_pe^{-(\alpha_ps)^2}\right)&=&-\frac{4}{3}\frac{\log^{3/2}(1/\epsilon_p)}{\alpha_p}(1+o(1)),
\end{eqnarray*}
where the \textcolor{black}{the first equality follows from \eqref{delta less than 1/2} and the second equation in \eqref{all def},} and the second equality follows from \eqref{log asymptotics} and the second equation in \eqref{all def}. 
The last term on the right hand side of \eqref{MWZ lower bound} is $o(1)$ because $\ell_p\tilde{\epsilon}_p=\frac{3\log^{1/2}(1/\epsilon_p)}{\alpha_p}\rightarrow\infty$. Therefore, overall we have the following lower bound on $\log \Pr(Z_1\leq\epsilon_p/2)$:
\begin{equation}
\label{final z1 tilde}
  \log\Pr(Z_1\leq\epsilon_p/2)=\log\Pr(\tilde{Z}_1\leq\tilde{\epsilon}_p)\geq -\frac{4}{3}\frac{\log^{3/2}(1/\epsilon_p)}{\alpha_p}(1+o(1))=-\frac{2}{3}\log (p) (1+o(1)).
\end{equation}

It remains to establish an upper bound on $\Pr(Z_2> \epsilon_p/2)$. 
Recall that a  centered random variable $X$ belongs to the sub-gamma
 family
${SG}(v,u)$ for $v, u > 0$ if
\begin{equation*}
    \log \mathbb{E} e^{tX} \leq \frac{t^2 v}{2(1 - tu)},
    \qquad
    \forall t: \; |t| < \frac{1}{u} \, .
  \end{equation*}
  If $X\in SG(v,u)$ then for every $t \geq 0$,
\begin{equation}
    \label{exponential ineq}
    \Pr(X>\sqrt{2vt}+ut)\leq e^{-t}.
\end{equation}
This inequality, as well as many other results related to sub-gamma random variables, can be found in chapter~{2.4} of \cite{boucheron2013concentration}.

  \textcolor{black}{As is well known, $C\times \chi^2(k)\in SG(2C^2k,2|C|)$, where $C$ is an arbitrary constant. Therefore, $Z_2-\mathbb{E}Z_2\in SG(\Delta^28(b_p+1),2\Delta)$, and since $\mathbb{E}Z_2=\Delta 4(b_p+1)$, we have}
%
for any $t\geq 0$, we have
\[
\Pr\left(Z_2>4\Delta(b_p+1)+2\Delta t+4\Delta\sqrt{(b_p+1)t}\right)\leq e^{-t}.
\]
Using \eqref{all def}, we obtain
\begin{equation}
\label{delta as power}
\Delta:=\frac{2^7\sigma_p^3}{np^2}e^{-\frac{p^2}{4\sigma_p^2}}<\mathrm{const}\times \frac{n^3}{\log^{3/2} p} p^{-\frac{p^2\pi^2}{n^22^5\beta^3}}\leq p^{-100}
\end{equation}
for some sufficiently small $\beta$ (that depends on $c$) and all sufficiently large $p$. 

Take $t=\epsilon_p/(5\Delta)$ and let $b_p=\lfloor \log p\rfloor$. \textcolor{black}{Then $t>p^\gamma$ for some $\gamma>0$ and all sufficiently large $p$, and hence $t\gg b_p$. Therefore, we have
\[
\epsilon_p/2=5\Delta t/2>2\Delta t+4\Delta(b_p+1)+4\Delta\sqrt{(b_p+1)t}
\]
for all sufficiently large $p$. This yields the following inequality}   
\begin{equation}
\label{final Z2}
\Pr(Z_2>\epsilon_p/2)< \Pr\left(Z_2>4\Delta(b_p+1)+2\Delta t+4\Delta\sqrt{(b_p+1)t}\right)\leq e^{-t}<e^{-p^\gamma}.
\end{equation}
Using \eqref{final z1 tilde} and \eqref{final Z2} in \eqref{split Y into Zs}, we obtain
\[
\Pr(Y_k\leq \epsilon_p)\geq p^{-\frac{2}{3}(1+o(1))}.
\]
\textcolor{black}{This implies that the term $\left(1-\Pr\left(Y_k\leq \epsilon_p\right)\right)^K$ on the right hand side of \eqref{the main inequality} converges to zero, because $K=\lfloor \frac{n}{4b_p+2}\rfloor-1\gg p^{\frac{2}{3}(1+o(1))}$.}
\vspace{2mm}

\textbf{Step 2: establishing an upper bound on $\Pr(R_k>\epsilon_p)$. } 
By definition of $R_k$, and recalling that \textcolor{black}{the indices are understood modulo $n$, we have
\begin{equation}
\label{Rk symmetric}
R_k=\sum_{s=b_p+1}^{n-b_p-1}2|y_{s+\tau_k}|^2d_{s},
\end{equation}
where $d_{s}:=\frac{1}{2n}|P_{\bu s}|^2$.  Suppose that $n$ is odd (the analysis for even $n$ is very similar and we omit it). Since  $\mathbb{E}|y_{s+\tau_k}|^2=1$ for any $s$ (see \eqref{y definition}), and $d_s=d_{n-s}$, we have
\begin{equation}
\label{E of R}
\mathbb{E}R_k=4\sum_{j=b_p+1}^{(n-1)/2}d_j.
\end{equation}
Further, from \eqref{y definition} and \eqref{Rk symmetric}, we have the representation
\begin{equation}
\label{Rk with gammas}
R_k=2\chi^2(1)d_{n-\tau_k}+\sum_{j=1}^{(n-1)/2}\gamma_j\chi^2_j(2)
\end{equation}
where $\gamma_j=\sum_{s\in S_j}d_s$ with $S_j$ being the set of all integer $s$ in between $b_p+1$ and $n-b_p-1$ such that $(s+\tau_k)\mod n$ equals $j$ or $n-j$. There is at most two elements in each $S_j$, $j=1,\dots,(n-1)/2$.}

\textcolor{black}{Now recall that if  $X\in SG\left( v_{X},u_{X} \right) $ and $Y\in SG\left( v_{Y},u_{Y}\right) $ are independent, then 
\begin{equation*}
X+Y\in SG\left( v_{X}+v_{Y},\max \left\{
u_{X},u_{Y}\right\} \right).
\end{equation*} 
Furthermore, if $X \in SG(v,u)$ then $X \in SG(v',u')$ for each $v' \geq v$, $u' \geq u$. 
These facts, representation \eqref{Rk with gammas}, and the inequality $\gamma_j^2\leq \sum_{s\in S_j}2d_{s}^2$ yield $R_k-\mathbb{E}R_k\in SG(v,u)$ with
\begin{eqnarray*}
v=\sum_{s=b_p+1}^{n-b_p-1}8d_s^2=\sum_{s=b_p+1}^{(n-1)/2}16d_s^2\qquad \text{and}\qquad u=\max_{b_p+1\leq s\leq n-b_p-1}4d_s=\max_{b_p+1\leq s\leq (n-1)/2}4d_s.
\end{eqnarray*}}

On the other hand, from \eqref{bound on d}, \eqref{all def} and \eqref{delta as power}, we have for all sufficiently large $p$,
\[
d_s\leq(2^4/\pi)\alpha_pe^{-\alpha_p^2s^2}+p^{-100}.
\]
This inequality together with \eqref{all def} and the fact that $b_p=\lfloor\log p\rfloor$ yield
\[
u<\frac{2^9\beta^{3/2}p^{-32\beta^{3}}}{\pi\log^{1/2}p}+4p^{-100}<p^{-\gamma}
\]
for some $\gamma>0$ (in the remaining part of the proof, $\gamma$ will denote a positive constant that may change its value from one appearance to another). Similarly, we obtain for all sufficiently large $p$,
\begin{eqnarray*}
v&<&16\int_{b_p}^{(n-1)/2}\left((2^4/\pi)\alpha_pe^{-\alpha_p^2x^2}+p^{-100}\right)^2\mathrm{d}x\\
&<&\gamma\alpha_p(1-\Phi(2\alpha_pb_p))+p^{-50}<p^{-\gamma},
\end{eqnarray*}
where $\Phi(\cdot)$ is the standard normal cdf. Finally, using $\eqref{E of R}$ we obtain
\[
\mathbb{E}R_k<4\int_{b_p}^{(n-1)/2}\left((2^4/\pi)\alpha_pe^{-\alpha_p^2x^2}+p^{-100}\right)\mathrm{d}x<p^{-\gamma}.
\]
Using these inequalities in \eqref{exponential ineq} for $X=R_k-\mathbb{E}R_k$, we obtain for any $t\geq 1$ and some $\gamma>0$,
\[
\Pr(R_k>p^{-\gamma}t)\leq e^{-t}.
\]
Taking $t=p^{\gamma}\epsilon_p$, and noting that $t>p^{\delta}$ for some small $\delta>0$ and all sufficiently large $p$, we conclude
\[
\Pr(R_k>\epsilon_p)\leq \exp\{-p^{\delta}\}.
\]
\textcolor{black}{This implies that the term $\sum_{k=1}^K\Pr\left( R_k> \epsilon_p\right)$ on the right hand side of \eqref{the main inequality} converges to zero, because $K=\lfloor \frac{n}{4b_p+2}\rfloor-1\ll \exp\{p^\delta\}$.}
\qed

\subsection{Proof details for Theorem \ref{thm: expectation}}\label{sec: thm2}


\textcolor{black}{As discussed in Section \ref{sec: main results}, it is sufficient to prove \eqref{to prove for thm2}.} Let $\eta_s = |y_s|^2 - \E |y_s|^2$. Then by \eqref{periodogam bound}, we have
\[
\lambda_p\leq \frac{1}{n}\sum_{s=0}^{n-1}|y_{s}|^2|P_{\bu s }|^2=\frac{1}{n}\sum_{s = 0}^{n - 1} |P_{\bu s}|^2 \E|y_s|^2 +\frac{1}{n} \sum_{s = 0}^{n - 1} |P_{\bu s}|^2 \eta_s =: c_n + Z_n,
\]
where $\bu$ is an arbitrary $p$-dimensional vector with $\|\bu\|=1$.
Since $\frac{1}{n}\sum_{s=0}^{n-1}|P_{\bu s}|^2=1$ and $\E|y_s|^2=1$, we have 
 $c_n\leq 1$. Therefore,
\bal{
\| \lambda_p \|_\kappa \leq 1 + \|Z_n\|_\kappa, \quad \kappa \geq 1,
}
 where $\|\xi\|_\kappa$ denotes the $L^\kappa$ norm of a random variable $\xi$, $\|\xi\|_\kappa = \big(\E |\xi|^k\big)^{1/\kappa}$.  It remains to bound $\|Z_n\|_\kappa$. Let the $\psi_1$-norm of a random variable $\xi$ be defined as 
\bal{
\| \xi \|_{\psi_1} = \inf\{t > 0: \E \exp (|\xi|/t) \leq 2\}.  
}
(The variables with finite $\psi_1$-norm are called sub-exponential.) Random variables $\eta_s$ are centered $\chi^2$ and therefore  sub-exponential,  and  $\| \eta_s\|_{\psi_1} \leq C < \infty$ for some absolute constant $C$.   By Theorem 2.8.2 in \cite{vershynin2018}, we have (assuming, for concreteness, that $n$ is odd)
\bal{
\Pr\left\{\Big|\sum_{s = 0}^{(n - 1)/2} w_s \eta_s\Big| \geq t\right\} \leq 
2 \exp \left( -\gamma \min \Big(\frac{t^2}{C^2 \|w\|_2^2}, \frac{t}{C \|w\|_\infty}\Big)\right),
}
where $\gamma>0$ is an absolute constant, $w_0=\frac{1}{n}|P_{\bu 0}|^2$ and $w_s=\frac{1}{n}\left(|P_{\bu s}|^2+|P_{\bu ,n-s}|^2\right)$ for $s>0$.

Let us take $\bu=p^{-1/2}(1,\dots,1)^\top$. Then it is straightforward to verify that
\[
\frac{1}{n}|P_{\bu 0}|^2=\frac{p}{n},\quad\text{and}\quad \frac{1}{n}|P_{\bu s}|^2=\frac{1}{pn}\left(\frac{\sin (\pi s p/n)}{\sin (\pi s/n)}\right)^2, 1\leq s\leq n-1.
\]
The latter expression is the value of the Fej\'{e}r kernel (divided ny $n$) at $2\pi s/n$. Since $\sin x\geq 2 x/\pi$ for all $x\in[0,\pi/2]$, we have 
\begin{equation*}
\frac{1}{n}|P_{\bu s}|^2\leq \frac{n}{4p}\frac{1}{s^2},\quad 1\leq s\leq (n-1)/2. 
\end{equation*}
Therefore, we see that $\max\{\|w\|_2^2, \|w\|_\infty\} \leq C_c$ for some constant $C_c$ that depends on $\lim p/n=c$. By Proposition 2.7.1 in \cite{vershynin2018}, this implies that $\big\|Z_n\big\|_{\psi_1} \leq C_c $, where the value of $C_c$ may change from one appearance to another, and 
\bal{
\Big\|Z_n\Big\|_\kappa \leq C \kappa \Big\|Z_n\Big\|_{\psi_1} \leq C_{\kappa,c},
}
which is what we wanted to prove. \qed

\subsection{Proof details for Theorem \ref{thm: lower bound}}\label{sec: thm3}

To complete the proof of Theorem \ref{thm: lower bound} outlined in Section \ref{sec: main results} we need, first, to establish \eqref{to show in sec 5}, and second, to show that the term $S_1$ in \eqref{referrred to next} is bounded away from zero.  

\textcolor{black}{Recall that \eqref{to show in sec 5} reads as 
$
\Pr\left(\min\{\xi_0,\xi_{(r)}\}>r/N\right)\rightarrow 1
$.
Since $\xi_0\sim \chi^2(1)$ and $r/N\rightarrow 0$, we have $\Pr\left(\xi_0>r/N\right)\rightarrow 1$. Further, as follows e.g.~from equation 1.9 of \cite{renyi53}, the order statistic $\xi_{(r)}$ has the following representation
\begin{equation}
\label{renyi}
\xi_{(r)}=\sum_{i=1}^{r}\frac{\zeta_i}{N-i},
\end{equation}
where $\zeta_{1},\dots,\zeta_{N-1}$ are mutually independent random variables with $\operatorname{Exp}(1/2)$ distribution. }
\textcolor{black}{ From \eqref{renyi}, we have \begin{equation*}
\mathbb{E}\xi_{(r)}=\sum_{i=1}^r\frac{2}{N-i},\qquad \operatorname{Var}\xi_{(r)}=\sum_{i=1}^{r}\frac{4}{(N-i)^2}.
\end{equation*}
In particular,
\begin{equation}
\label{mean var}
\mathbb{E}\xi_{(r)}>2r/N,\qquad \operatorname{Var}\xi_{(r)}<5r/N^2
\end{equation}
for all sufficiently large $p$ (recall that $r:=\lfloor p^{1-1/m-\epsilon/2}\rfloor$ and $N=(n+1)/2$ with $p/n\rightarrow c\in(0,1]$). Therefore, by Chebyshev's inequality,
\begin{eqnarray*}
    \Pr\left(\xi_{(r)}<r/N\right)<\Pr\left(\xi_{(r)}-\mathbb{E}\xi_{(r)}<-r/N\right)\leq\frac{5r/N^2}{r^2/N^2}\rightarrow 0.
\end{eqnarray*}
Hence $\Pr(\xi_{(r)}>r/N)\rightarrow 1$, which completes the proof of \eqref{to show in sec 5}. 
}

We now turn to the proof of the boundedness of $S_1$ away from zero. Recall that 
\[
S_1:=\frac{1}{2}T_{p-1}(0)+\sum_{j=r+1}^{N-1}T_{p-1}(\omega_{\sigma(j)}),\qquad\text{and let }\qquad S_2:=\sum_{j=1}^{r}T_{p-1}(\omega_{\sigma(j)}).
\]
As follows from \eqref{discrete trigonometric normalization}, $S_2=1/2-S_1$.  We show that $S_1$ is bounded away from zero by using the $L_1$ version of the Bernstein inequality for the derivative of a cosine trigonometric polynomial (see e.g.~Theorem 2.5 in Chapter 4 of \cite{devore93}): 
\begin{equation}
\label{bernstein intro}
\int_0^{\pi}\left|T_{p-1}'(\omega)\right|\mathrm{d}\omega\leq (p-1)\int_0^{\pi}|T_{p-1}(\omega)|\mathrm{d}\omega=(p-1)\frac{\pi}{n}.
\end{equation}
Here the last equality follows from \eqref{trigonometric normalization}.

Let us call $\omega_i$ type 2 if $i=\sigma(j)$ with $1\leq j\leq r$ and type 1 otherwise. Then the set $\{\omega_0,\omega_1,\dots,\omega_{N-1}\}$ splits into a union of  $K$ interlacing clusters (a cluster may contain only one element) of type 1 and type 2. The first cluster is of type 1 because $\omega_0$ is of type 1. Assuming that  none of $\sigma(j):1\leq j\leq r$ equals $N-1$ on $A_p$ (such an assumption is without loss of generality because the probability of the latter event equals $(N-r-1)/(N-1)$, which converges to 1), the last cluster is also of type 1 because $\omega_{N-1}$ is of type 1. In particular, $K$ is an odd integer.

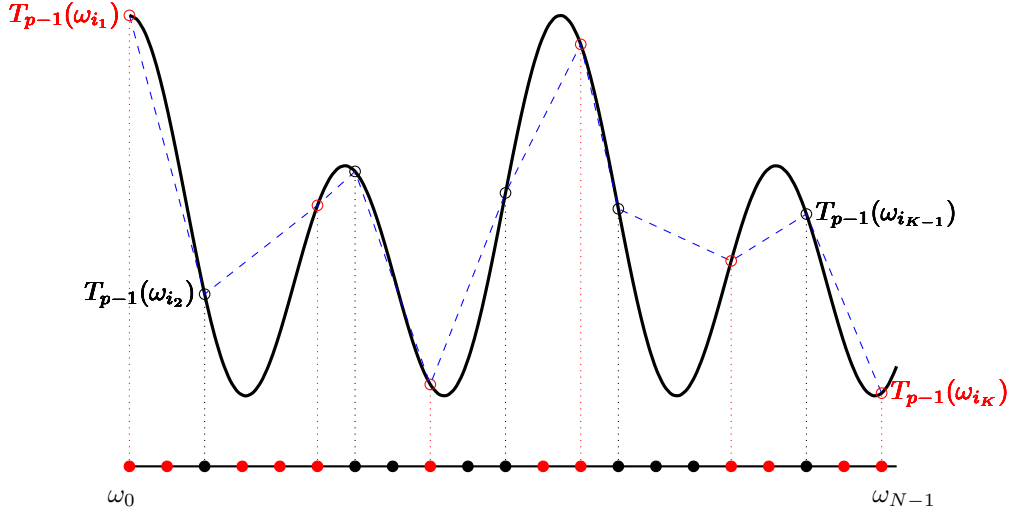
\begin{figure}
    \centering
    \begin{tikzpicture}[domain=0:10.2]
    \draw[thick, -] (0,0) -- (10.2,0);
    \foreach \x in {0,1,3,4,5,8,11,12, 16,17, 19,20} \filldraw[red] ({\x/2},0) circle (2pt);
    \foreach \x in {2,6,7,9,10,13,14,15,18} \filldraw ({\x/2},0) circle (2pt);
    \draw (-0.1,-0.2) node[anchor=north]{$\omega_0$};
    \draw (10.3,-0.2) node[anchor=north]{$\omega_{N-1}$};
    \draw[very thick, samples=100] plot [smooth] (\x,{3+cos(20*pi*\x)+2*cos(40*pi*\x)});
    \foreach \x in {0,2.5,4,6,8,10} \draw[dotted,red] (\x,0)--(\x,{3+cos(20*pi*\x)+2*cos(40*pi*\x)});
    \foreach \x in {1,3,5,6.5,9} \draw[dotted] (\x,0)--(\x,{3+cos(20*pi*\x)+2*cos(40*pi*\x)});
    \foreach \x [count=\n] in {0,1,2.5,3,4,5,6,6.5,8,9,10} {
        \pgfmathsetmacro{\y}{3+cos(20*pi*\x)+2*cos(40*pi*\x)} 
        \ifnum\n=1 
            \xdef\plotcoordinates{(\x,\y)}
        \else 
            \xdef\plotcoordinates{\plotcoordinates (\x,\y)}
        \fi;
        \draw[red] (0,{3+cos(20*pi*0)+2*cos(40*pi*0)}) node[anchor=east]{$T_{p-1}(\omega_{i_1})$};
        \draw (1,{3+cos(20*pi*1)+2*cos(40*pi*1)}) node[anchor=east]{$T_{p-1}(\omega_{i_2})$};
        \draw[red] (10,{3+cos(20*pi*10)+2*cos(40*pi*10)}) node[anchor=west]{$T_{p-1}(\omega_{i_K})$};
        \draw (9,{3+cos(20*pi*9)+2*cos(40*pi*9)}) node[anchor=west]{$T_{p-1}(\omega_{i_{K-1}})$};
    }

    \draw[dashed, blue] plot coordinates {\plotcoordinates};
    \foreach \x in {0,2.5,4,6,8,10} {
        \pgfmathsetmacro{\y}{3+cos(20*pi*\x)+2*cos(40*pi*\x)} 
        \draw[red] (\x,\y) circle (2pt);}
        \foreach \x in {1,3,5,6.5,9} {
        \pgfmathsetmacro{\y}{3+cos(20*pi*\x)+2*cos(40*pi*\x)} 
        \draw (\x,\y) circle (2pt);}
    
\end{tikzpicture}
    \caption{Illustration of the inequality \eqref{trigonometric variability}. Adjacent red dots form clusters of $\omega_i$ of the first type. Adjacent black dots form clusters of $\omega_i$ of the second type. Total variation of the smooth graph, representing the left hand of \eqref{trigonometric variability}, is larger than that of the piece-wise linear dashed graph, representing  the right hand side of \eqref{trigonometric variability}. }
    \label{fig:illustr}
\end{figure}

Picking just one element from each of the clusters, we obtain a sequence $\omega_{i_1},\omega_{i_2},\dots,\omega_{i_{K}}$ with odd elements of type 1 and even elements of type 2. Clearly, 
\begin{eqnarray}
\label{trigonometric variability}
\int_0^{\pi}\left|T_{p-1}'(\omega)\right|\mathrm{d}\omega&\geq& |T_{p-1}(\omega_{i_1})-T_{p-1}(\omega_{i_2})|+\dots+|T_{p-1}(\omega_{i_{K-1}})-T_{p-1}(\omega_{i_{K}})|.
\end{eqnarray}
Figure \ref{fig:illustr} illustrates this inequality. On the other hand, the right hand side of \eqref{trigonometric variability} is no smaller than
\begin{eqnarray*}
    &&\sum_{\text{even }0<k<K}(T_{p-1}(\omega_{i_k})-T_{p-1}(\omega_{i_{k+1}}))-\sum_{\text{odd }0<k<K}(T_{p-1}(\omega_{i_k})-T_{p-1}(\omega_{i_{k+1}}))\\
    &\geq&2\sum_{\text{even }0<k< K}T_{p-1}(\omega_{i_k})-\left(T_{p-1}(\omega_{i_1})+2\sum_{\text{odd }1<k\leq K}T_{p-1}(\omega_{i_k})\right).
\end{eqnarray*}
Since all $\omega_{i_k}$ with odd $k$ are of type 1 and $T_{p-1}$ is a non-negative trigonometric polynomial, the expression in the large brackets in the latter display is no larger than $2S_1$. Hence, combining the latter inequality with \eqref{bernstein intro} and \eqref{trigonometric variability} and recalling that all $\omega_{i_k}$ with even $k$ are of type 2, we obtain
\begin{equation}
\label{bound on linear comb of S1 and S2}
\frac{(p-1) \pi}{n}\geq 2\sum_{\omega_{i_k}\text{ of type 2}}T_{p-1}(\omega_{i_k})-2S_1.
\end{equation}

Let us now show that  $\omega_{i}$ of type 2 cannot cluster in groups of $m$ or more adjacent $\omega_i$'s with probability approaching one as $p\rightarrow\infty$. Recall that $\omega_{i}$ is of type 2 if $i=\sigma(j)$ with $1\leq j\leq r$, and that $\sigma$ is a random permutation uniformly distributed on the set of all permutations of the integers $1,\dots,N-1$. Consider  the probability that randomly picked  $r$ integers out of $1,\dots, N-1$ contain at least $m$ adjacent integers, that is, integers $i+j$ with $j=1,\dots,m$ for some $0\leq i\leq N-1-m$. We denote this probability as $P(N,r,m)$. We need to prove that $P(N,r,m)\rightarrow 0$ as $N\rightarrow \infty$.

Denote the number of possible ordered sequences of $r$ integers out of $N-1$ that would contain at least $m$ adjacent integers as $\gamma(N,r,m)$. Then,
\[
P(N,r,m)=\frac{\gamma(N,r,m)}{(N-1)!/(N-1-r)!}.
\]
On the other hand, $\gamma(N,r,m)$ is no larger than the product of the following three numbers: $r!/(r-m)!$ - the number of possible choices of $m$ ordered places out of $r$; $N-m$ - the number of possible choices of $m$ adjacent integers out of $1,\dots, N-1$; and $(N-1-m)!/(N-1-r)!$ - the number of possible choices of ordered sequences of length $r-m$ out of remaining $N-1-m$ integers. Therefore,
\[
P(N,r,m)\leq \frac{\frac{r!}{(r-m)!}\times (N-m)\times \frac{(N-1-m)!}{(N-1-r)!}}{\frac{(N-1)!}{(N-1-r)!}}=\frac{r\times\dots\times(r-m+1)}{(N-1)\times\dots\times(N-m+1)}.
\]
Clearly, as $N\rightarrow\infty$ while $m$ is being fixed, the right hand side converges to zero as long as $r\ll N^{\frac{m-1}{m}}$, which is indeed the case for  $r:=\lfloor p^{1-1/m-\epsilon/2}\rfloor$.

To summarize, as we have just shown, $\omega_{i}$ of type 2 cannot cluster in groups of $m$ or more adjacent $\omega_i$'s with probability approaching one as $p\rightarrow\infty$. Therefore, returning to inequality \eqref{bound on linear comb of S1 and S2}, if we pick $\omega_{i_k}$ of the second types so that the corresponding $T_{p-1}(\omega_{i_k})$ is the maximum possible in the cluster that $\omega_{i_k}$ belongs to (as illustrated in Figure \ref{fig:illustr}), we have, on a sequence of high probability events, 
\begin{equation*}
\sum_{\omega_{i_k}\text{ of type 2}}T_{p-1}(\omega_{i_k})\geq \frac{1}{m-1}S_2.
\end{equation*}
Combining this with \eqref{bound on linear comb of S1 and S2} and using the fact that $S_2=1/2-S_1$, we obtain 
\[
S_1\geq \frac{m-1}{2m}\left(\frac{1}{m-1}-\frac{(p-1)\pi}{n}\right).
\]
But the right hand side is bounded away from zero for all sufficiently large $p$ as long as $\lim n/p>\pi(m-1)$, which is required by Theorem \ref{thm: lower bound}. This concludes the proof.
\qed

\subsection{Proof details for Theorem \ref{thm: small T/p}}\label{sec: thm4}

We need the following result (see its proof in the Appendix), showing that trigonometric polynomials cannot decay too fast  near the point of their maximum. 
\begin{lemma}
\label{lem: Stechkin 1}
    Let
$T_k(\omega)=a_0+\sum_{j=1}^k a_j\cos(j\omega)$ 
be a non-negative trigonometric cosine polynomial of degree $k$ with real coefficients, and suppose that $\max T_k(\omega)=:M>0$. If $T_k(\bar{\omega})=M$, then \[
T_k(\bar{\omega}+\eta)\geq M\cos^2\frac{k\eta}{2}\qquad \text{for all }|\eta|\leq \pi/k.\]
\end{lemma}

The setting of our
proof of Theorem \ref{thm: small T/p} is exactly the same as that of the proof of Theorem \ref{thm: lower bound} above. In particular, we keep the definitions of $r$, of type 1 and type 2 clusters, and of the events $A_p$. Now let $r_1:=\lfloor p^{1/2-1/(2m)-\epsilon/4}\rfloor$ so that $r_1/\sqrt{r}\rightarrow 1$, and assume that events $A_p$ are adjusted to imply that $\xi_{(r_1)}>r_1/N$. Such an adjustment is feasible because $\Pr(\xi_{(r_1)}>r_1/N)\rightarrow 1$, as can be shown similarly to  \eqref{to show in sec 5} above.
\textcolor{black}{Let 
\[
S_{11}:=\frac{1}{2}T_{p-1}(0)+\sum_{j=r_1+1}^{N-1}T_{p-1}(\omega_{\sigma(j)}),\qquad S_{21}:=\sum_{j=1}^{r_1}T_{p-1}(\omega_{\sigma(j)})=1/2-S_{11},
\]
where the latter equality follows from \eqref{discrete trigonometric normalization}. Similarly to \eqref{referrred to next}, on $A_p$ we have
\begin{equation}
\label{first component of max}
\lambda_p>\frac{r_1}{N}S_{11}.
\end{equation}}

Further, let $M:=\max_\omega T_{p-1}(\omega)=T_{p-1}(\bar{\omega})$. Since $T_{p-1}(\omega)$ is an even function, we can assume without loss of generality that $\bar{\omega}\in[0,\pi]$. Clearly,
$M\geq S_{21}/r_1=(1/2-S_{11})/r_1$.
This inequality and Lemma \ref{lem: Stechkin 1} imply that
\begin{equation}
\label{implication of Lemma1}
T_{p-1}(\bar{\omega}+\eta)\geq \frac{1/2-S_{11}}{r_1}\cos^2\frac{p\eta}{2}\qquad\text{for all }|\eta|\leq\frac{\pi}{p}.
\end{equation}
(Strictly following Lemma \ref{lem: Stechkin 1} yields even a stronger statement with $\cos^2\frac{(p-1)\eta}{2}$ and $|\eta|\leq\frac{\pi}{p-1}$ in the above display.)

Consider the following segment with a center at $\bar{\omega}$: 
\[
I:= \left[\bar{\omega}-(1-\delta_0)\pi/p,\bar{\omega}+(1-\delta_0)\pi/p\right],
\]
where constant $\delta_0>0$ is such that, for all sufficiently large $p$, $n/p\geq m/(1-\delta_0)$ (the existence of such a constant follows from the conditions of Theorem \ref{thm: small T/p}).  Inequality \eqref{implication of Lemma1} implies that, for any point $\omega_{\sigma(j)}:=\frac{2\pi \sigma(j)}{n}\in I$, we must have 
\[
T_{p-1}(\omega_{\sigma(j)})\geq \frac{1/2-S_{11}}{r_1}\cos^2\frac{(1-\delta_0)\pi}{2}.
\]
Note that the number of adjacent points of form $2\pi s/n$, $s\in\mathbb{Z}$ in $I$ cannot be smaller than
\[
\left\lfloor\frac{(1-\delta_0)2\pi/p}{2\pi/n}\right\rfloor=\left\lfloor\frac{(1-\delta_0)n}{p}\right\rfloor\geq m.
\]
Therefore, if $I\subseteq [\omega_1,\omega_{N-1}]$, at least one of the points $\omega_s\in I$ must be of type 1, and hence, by \eqref{lambdap for thm3} we must have
\[
\lambda_p\geq \xi_{(r)}T_{p-1}(\omega_s)> \frac{r}{N}\frac{1/2-S_{11}}{r_1}\cos^2\frac{(1-\delta_0)\pi}{2}.
\]
Combining this with \eqref{first component of max}, we obtain
\[
\lambda_p>\max\left\{\frac{r_1}{N}S_{11},\frac{r}{N}\frac{1/2-S_{11}}{r_1}\cos^2\frac{(1-\delta_0)\pi}{2}\right\}.
\]
Since $S_{11}\in[0,1/2]$, we conclude that $\lambda_p>C r_1/N$ for some positive constant $C$, which implies that, on $A_p$, $\lambda_p>p^{-1/2-1/(2m)-\epsilon}$ for all sufficiently large $p$, and Theorem \ref{thm: small T/p} holds.

To conclude the proof, we need to consider cases when $I\nsubseteq [\omega_1,\omega_{N-1}]$. This is possible only if $\bar{\omega}\in[0,(1-\delta_0)\pi/p]$ or if $\bar{\omega}\in [\pi-(1-\delta_0)\pi/p,\pi]$. Obviously, in the former case, $\omega_1\in I$, whereas in the latter case, $\omega_{N-1}\in I$. Let us adjust events $A_p$ so that they imply that both $\omega_1$ and $\omega_{N-1}$ are of type 1. Then, we can proceed as in the case $I\subseteq [\omega_1,\omega_{N-1}]$, setting $\omega_s$ to $\omega_1$ or $\omega_{N-1}$.  \qed

\section{Conclusion}\label{sec: conclusion}
This paper proves that, somewhat unexpectedly, the smallest non-zero eigenvalues of the rectangular $p\times n$ Toepliz and circulant matrices, \eqref{toeplitz matrix} and \eqref{circulant matrix}, converge in probability and in expectation to zero as $p\rightarrow \infty$ for all $c:=\lim p/n\in(0,1]$. We show that the rate of the convergence is faster than poly-log, and derive some polynomial upper bounds on the rate. Our Monte Carlo exercises suggest that the actual rate might be polynomial, and that the upper bounds provided by Theorems \ref{thm: lower bound} and \ref{thm: small T/p} could be sub-optimal.

In future research, it would be valuable to prove that the rate of convergence is indeed polynomial and characterize it precisely. Additionally, it would be interesting to extend our results to the case of non-Gaussian entries of the Toeplitz and circulant matrices.

\begin{appendix}
\section*{}
\subsection*{Proof of Lemma \ref{lem: barnett}}
Recall the Poisson summation formula:
\begin{equation}
\label{Poisson summation}
\sum_{j\in\mathbb{Z}}e^{-2\pi\rmi j\omega}f(j)=\sum_{m\in\mathbb{Z}}\hat{f}(\omega+m),
\end{equation}
where $f(x)$ is any integrable function and $\hat{f}(\omega)=\int f(x)e^{-2\pi\rmi\omega x}\mathrm{d}x$ is its Fourier transform.
    Applying \eqref{Poisson summation} to the Fourier pair
    \[
    f(x)=(2\pi\sigma_p^2)^{-1/2}\exp\left\{-\frac{(x-\lceil p/2\rceil)^2}{2\sigma_p^2}\right\}\quad\text{and}\quad \hat{f}(\omega)=\exp\left\{-2\pi\rmi \lceil p/2\rceil\omega-2(\pi\sigma_p\omega)^2\right\},
    \]
    we obtain
    \[
    (2\pi\sigma_p^2)^{-1/2}\sum_{j\in\mathbb{Z}}e^{-2\pi\rmi j\omega}\exp\left\{-\frac{(j-\lceil p/2\rceil)^2}{2\sigma_p^2}\right\}=e^{-2\pi\rmi \lceil p/2\rceil \omega}\sum_{m\in\mathbb{Z}}\exp\left\{-2(\pi\sigma_p(\omega+m))^2\right\}.
    \]
    Setting $\omega=-k/n$ and grouping terms on the left hand side yields \eqref{branett transformed}.

\subsection*{Proof of Lemma \ref{lem: exponential d}}
    First, borrowing some ideas from the proof of Proposition 7 of \cite{barnett22}, we bound $|\hat{w}_s|$, defined in \eqref{branett transformed}. We have
    \[
    |\hat{w}_s|=\sum_{m\in\mathbb{Z}}e^{-2\left(\frac{\pi\sigma_p}{n}\right)^2(s+mn)^2}=\sum_{m\geq 1}e^{-2\left(\frac{\pi\sigma_p}{n}\right)^2(s+mn)^2}+e^{-2\left(\frac{\pi\sigma_p}{n}\right)^2s^2}+\sum_{m\leq -1}e^{-2\left(\frac{\pi\sigma_p}{n}\right)^2(s+mn)^2}.
    \]
    For the first sum on the right hand side, we have
    \begin{eqnarray}\notag
    \sum_{m\geq 1}e^{-2\left(\frac{\pi\sigma_p}{n}\right)^2(s+mn)^2}&<&\sum_{m\geq 1}e^{-2\left(\frac{\pi\sigma_p}{n}\right)^2(s^2+2sn+(mn)^2)}\\\notag
    &=&e^{-2\left(\frac{\pi\sigma_p}{n}\right)^2(s^2+2sn)}\sum_{m\geq 1}e^{-2\left(\frac{\pi\sigma_p}{n}\right)^2(mn)^2}\\\notag
    &<&e^{-2\left(\frac{\pi\sigma_p}{n}\right)^2(s^2+2sn)}\int_0^\infty e^{-2\pi^2\sigma_p^2x^2}\mathrm{d}x\\\label{argument}
    &=&\frac{1}{\sqrt{8\pi}\sigma_p}e^{-2\left(\frac{\pi\sigma_p}{n}\right)^2(s^2+2sn)}.
    \end{eqnarray}
    If $0\leq s\leq n/2$, then 
    \begin{eqnarray*}
    \sum_{m\leq -1}e^{-2\left(\frac{\pi\sigma_p}{n}\right)^2(s+mn)^2}&=&\sum_{m\geq 0}e^{-2\left(\frac{\pi\sigma_p}{n}\right)^2(n-s+mn)^2}\\
    &\leq&e^{-2\left(\frac{\pi\sigma_p}{n}\right)^2s^2}+\sum_{m\geq 1}e^{-2\left(\frac{\pi\sigma_p}{n}\right)^2(n-s+mn)^2}.
    \end{eqnarray*}
    If $0\leq s\leq n /2$, the elements of the latter sum are no larger than those of\linebreak $\sum_{m\geq 1}e^{-2\left(\frac{\pi\sigma_p}{n }\right)^2(s+mn )^2}$. Hence, the argument that led us to \eqref{argument} yields
    \[
    \sum_{m\leq -1}e^{-2\left(\frac{\pi\sigma_p}{n }\right)^2(s+mn )^2}\leq e^{-2\left(\frac{\pi\sigma_p}{n }\right)^2s^2}+\frac{1}{\sqrt{8\pi}\sigma_p}e^{-2\left(\frac{\pi\sigma_p}{n }\right)^2(s^2+2sn )}.
    \]
    Summing up, for $0\leq s\leq n /2$ we have
    \begin{equation}
    \label{bound on vkhat}
    |\hat{w}_s|\leq 2e^{-2\left(\frac{\pi\sigma_p}{n }\right)^2s^2}\left(1+\frac{e^{-4s\pi^2\sigma_p^2/n }}{\sqrt{8\pi}\sigma_p}\right).
    \end{equation}
    By definition, it is clear that, for $n /2<s\leq n -1$, $|\hat{w}_s|=|\hat{w}_{n -s}|$.

    We can similarly obtain an upper bound on $|w_j|$, defined in \eqref{barnett direct}. For $\lceil p/2\rceil\leq j\leq \lceil p/2\rceil +n /2$, we have
    \begin{eqnarray*}
    \sum_{s\geq 0}e^{-\frac{1}{2}(j-\lceil p/2\rceil+sn )^2/\sigma_p^2}&\leq& e^{-\frac{1}{2}(j-\lceil p/2\rceil)^2/\sigma_p^2}\left(1+\sum_{s\geq 1}e^{-\frac{1}{2}(sn )^2/\sigma_p^2}\right)\\
    &\leq&e^{-\frac{1}{2}(j-\lceil p/2\rceil)^2/\sigma_p^2}\left(1+\int_0^\infty e^{-\frac{1}{2}(xn )^2/\sigma_p^2}\mathrm{d}x\right)\\
    &=&e^{-\frac{1}{2}(j-\lceil p/2\rceil)^2/\sigma_p^2}\left(1+\frac{\sqrt{2\pi}\sigma_p}{2n }\right).
    \end{eqnarray*}
    Further,
    \begin{eqnarray*}
        \sum_{s< 0}e^{-\frac{1}{2}(j-\lceil p/2\rceil+sn )^2/\sigma_p^2}&=&\sum_{s\geq 0}e^{-\frac{1}{2}(n +\lceil p/2\rceil-j+sn )^2/\sigma_p^2}.
    \end{eqnarray*}
    Since by assumption, $n +\lceil p/2\rceil -j\geq j-\lceil p/2\rceil\geq 0$, the latter sum is no larger than the sum studied in the preceding display. Hence overall
    \[
    \sum_{s\in\mathbb{Z}}e^{-\frac{1}{2}(j-\lceil p/2\rceil+sn )^2/\sigma_p^2}\leq 2e^{-\frac{1}{2}(j-\lceil p/2\rceil)^2/\sigma_p^2}\left(1+\frac{\sqrt{2\pi}\sigma_p}{2n }\right),
    \]
    and thus, for $\lceil p/2\rceil\leq j\leq \lceil p/2\rceil +n /2$, we have
    \begin{equation}
    \label{bound on vj}
    |w_j|\leq (2\pi\sigma_p^2)^{-1/2}e^{-\frac{1}{2}(j-\lceil p/2\rceil)^2/\sigma_p^2}\left(2+\frac{\sqrt{2\pi}\sigma_p}{n }\right).
    \end{equation}

    Now note that for any $s=0,\dots, n-1$,
    \[
    \left|\sum_{j=p}^{n-1}e^{\frac{2\pi\rmi}{n}sj}w_j\right|\leq \sum_{j=p}^{n-1}|w_j|\leq 2\sum_{j=p}^{\lfloor\lceil p/2\rceil+n/2\rfloor}|w_j|.
    \]
    The latter inequality follows from the fact that, by definition \eqref{barnett direct}, $w_j=w_{2\lceil p/2\rceil+n-j}$. In particular, for any $j$ such that $\lfloor\lceil p/2\rceil+n/2\rfloor<j\leq n-1$, we have $w_j=w_k$ with $k:=2\lceil p/2\rceil +n-j$ satisfying $p\leq k\leq \lfloor\lceil p/2\rceil+n/2\rfloor$.

    Using bound \eqref{bound on vj} in the above display, we obtain
    \begin{eqnarray*}
    \left|\sum_{j=p}^{n-1}e^{\frac{2\pi\rmi}{n}sj}w_j\right|&<&2(2\pi\sigma_p^2)^{-1/2}\left(2+\frac{\sqrt{2\pi}\sigma_p}{n}\right)\sum_{j=p}^\infty e^{-\frac{1}{2}(j-\lceil p/2\rceil)^2/\sigma_p^2}\\
    &<&2(2\pi\sigma_p^2)^{-1/2}\left(2+\frac{\sqrt{2\pi}\sigma_p}{n}\right)\int_{p/2-2}^\infty e^{-\frac{1}{2}x^2/\sigma_p^2}\mathrm{d}x.
    \end{eqnarray*}
    This yields
    \[
    \left|\sum_{j=p}^{n -1}e^{\frac{2\pi\rmi}{n }sj}w_j\right|\leq\left(4+\frac{\sqrt{8\pi}\sigma_p}{n }\right)\left(1-\Phi\left(\frac{p-4}{2\sigma_p}\right)\right),
    \]
    where $\Phi$ denotes the cumulative distribution function of the standard normal random variable.
    Since for any $t>0$, 
    \begin{equation}
    \label{elementary ineq}
    1-\Phi(t)<\frac{1}{\sqrt{2\pi}}\frac{1}{t}e^{-t^2/2},
    \end{equation}
    we conclude that for $p>4$,
    \begin{equation}
        \label{bound on the rest of DFT}
        \left|\sum_{j=p}^{n-1}e^{\frac{2\pi\rmi}{n}sj}w_j\right|\leq\left(\frac{8}{\sqrt{2\pi}}+\frac{4\sigma_p}{n}\right)\frac{\sigma_p}{p-4}e^{-\frac{(p-4)^2}{8\sigma_p^2}}.
    \end{equation}

    Next, let us establish a lower bound on $\sqrt{\sum_{j=0}^{p-1}w_j^2}$. From the definition of $w_j$ we have
    \[
    w_j>(2\pi\sigma_p^2)^{-1/2}e^{-\frac{1}{2}(j-\lceil p/2\rceil)^2/\sigma^2_p}.
    \]
    Therefore,
    \begin{eqnarray*}
    \sum_{j=0}^{p-1}w_j^2&>&(2\pi\sigma_p^2)^{-1}\sum_{j=0}^{p-1}e^{-(j-\lceil p/2\rceil)^2/\sigma^2_p}\\
    &>&(2\pi\sigma_p^2)^{-1}\left(\int_{-p/2}^{p/2}e^{-x^2/\sigma_p^2}\mathrm{d}x-2\right)\\
    &=&\frac{1}{\sqrt{\pi}\sigma_p}\Phi\left(\frac{p}{\sqrt{2}\sigma_p}\right)-\frac{1}{2\sqrt{\pi}\sigma_p}-\frac{1}{\pi\sigma_p^2}.
    \end{eqnarray*}
    Using \eqref{elementary ineq}, we arrive at the inequality
    \begin{equation}
    \label{lower bound on vbar}
    \sum_{j=0}^{p-1}w_j^2>\frac{1}{2\sqrt{\pi}\sigma_p}-\frac{1}{\pi\sigma_p^2}-\frac{1}{\pi p}e^{-p^2/(4\sigma_p^2)}.
    \end{equation}

    By definition, 
    \begin{eqnarray*}
    \frac{1}{n}|P_{\bu s}|^2=\frac{1}{n}\left|(\sqrt{n}\mathcal{F}\bu_C)_s\right|^2&=& \frac{1}{n\sum_{j=0}^{p-1}w_j^2}\left|\hat{w}_s-\sum_{j=p}^{n-1}e^{\frac{2\pi\rmi}{n}sj}w_j\right|^2.
    \end{eqnarray*}
    Combining \eqref{bound on vkhat}, \eqref{bound on the rest of DFT} and \eqref{lower bound on vbar}, we arrive at the following inequality for any integer $s\in[0,n/2]$ and $p>4$
    \[
    \frac{1}{n}|P_{\bu s}|^2<\frac{8e^{-4\left(\frac{\pi\sigma_p}{n }\right)^2s^2}\left(1+\frac{e^{-4s\pi^2\sigma_p^2/n }}{\sqrt{8\pi}\sigma_p}\right)^2+2\left(\frac{8}{\sqrt{2\pi}}+\frac{4\sigma_p}{n}\right)^2\frac{\sigma_p^2}{(p-4)^2}e^{-\frac{(p-4)^2}{4\sigma_p^2}}}{n\left(\frac{1}{2\sqrt{\pi}\sigma_p}-\frac{1}{\pi\sigma_p^2}-\frac{1}{\pi p}e^{-p^2/(4\sigma_p^2)}\right)}.
    \]
    Since $p,n\rightarrow_c\infty$ and $\sigma_p\rightarrow\infty$ so that $\sigma_p=o(p)$, we have for all sufficiently large $p$ and any $s\in[0,n/2]$
    \begin{eqnarray*}
        \left(1+\frac{e^{-4s\pi^2\sigma_p^2/n }}{\sqrt{8\pi}\sigma_p}\right)^2&\leq&\sqrt{2},\\
        \left(\frac{8}{\sqrt{2\pi}}+\frac{4\sigma_p}{n}\right)^2\frac{\sigma_p^2}{(p-4)^2}&\leq&\sqrt{2}\frac{64}{2\pi}\frac{\sigma_p^2}{p^2},\\
        \frac{1}{2\sqrt{\pi}\sigma_p}-\frac{1}{\pi\sigma_p^2}-\frac{1}{\pi p}e^{-p^2/(4\sigma_p^2)}&\geq& \frac{1}{\sqrt{2}}\frac{1}{2\sqrt{\pi}\sigma_p},
    \end{eqnarray*}
    which (together with $1<\sqrt{\pi}<2$ and $p=o(\sigma_p^2)$) yields 
    \[
    \frac{1}{n}|P_{\bu s}|^2<\frac{2^6\sigma_p}{n}e^{-4\left(\frac{\pi\sigma_p}{n }\right)^2s^2}+\frac{2^8 \sigma_p^3}{n p^2}e^{-\frac{p^2}{4\sigma_p^2}}.
    \]
    Finally, the equality $|P_{\bu s}|=|P_{\bu,n-s}|=|P_{\bu,-s}|$ follows directly from the definition of $\mb{P}_\bu$ as the DFT of $\bu_C$.
  
    \subsection*{Proof of Lemma \ref{lem: prelim}} 
For any positive integer $J$, we have
\begin{eqnarray*}
\sum_{j=0}^\infty\log\left(1+\ell_p \alpha_pe^{-(\alpha_p j)^2}\right)&\geq&\sum_{j=0}^{J} \left(\log (\ell_p \alpha_p) -(\alpha_pj)^2\right)\\&=&(J+1)\log (\ell_p \alpha_p) -\frac{\alpha_p^2}{6}J(J+1)(2J+1).
\end{eqnarray*}
Setting
$
J=\left\lfloor \frac{\log^{1/2}(\ell_p\alpha_p)}{\alpha_p}\right\rfloor$,
we obtain
\begin{equation}
\label{log lower bound}
\sum_{j=0}^\infty\log\left(1+\ell_p \alpha_pe^{-(\alpha_p j)^2}\right)\geq \frac{2}{3}\frac{\log^{3/2}(\ell_p\alpha_p)}{\alpha_p}(1+o(1)).
\end{equation}
On the other hand, for the above choice of $J$, we have $\ell_p\alpha_pe^{-(\alpha_p J)^2}\geq 1$ and therefore,
\begin{eqnarray*}
\sum_{j=0}^\infty\log\left(1+\ell_p \alpha_pe^{-(\alpha_p j)^2}\right)&\leq& \sum_{j=0}^{J} \log \left(2\ell_p \alpha_pe^{-(\alpha_pj)^2}\right)+\int_{J}^\infty\log\left(1+\ell_p \alpha_p e^{-(\alpha_p x)^2}\right)\mathrm{d}x\\
&=&\frac{2}{3}\frac{\log^{3/2}(\ell_p\alpha_p)}{\alpha_p}(1+o(1))+\int_{J}^\infty\log\left(1+\ell_p \alpha_p e^{-(\alpha_p x)^2}\right)\mathrm{d}x.
\end{eqnarray*}
Let us show that the latter integral is of order lower than $\log^{3/2}(\ell_p\alpha_p)/\alpha_p$. Integrating by parts, we obtain
\begin{eqnarray}
\notag
&&\int_{J}^\infty\log\left(1+\ell_p \alpha_p e^{-(\alpha_p x)^2}\right)\mathrm{d}x\\\notag &=&-J\log(1+\ell_p \alpha_p e^{-(\alpha_pJ)^2})+2\alpha_p^3\ell_p\int_{J}^\infty\frac{ x^2e^{-(\alpha_px)^2}}{1+\ell_p \alpha_p e^{-(\alpha_p x)^2}}\mathrm{d}x\\\notag
&\leq&2\alpha_p^3\ell_p\int_{J}^\infty x^2e^{-(\alpha_px)^2}\mathrm{d}x=\ell_p \Gamma\left(3/2,(\alpha_pJ)^2\right)\\\notag &=&J\ell_p\alpha_p e^{-(\alpha_pJ)^2}(1+o(1)).
\end{eqnarray}
Here $\Gamma(a,z):=\int_z^\infty e^{-t}t^{a-1}dt$ is the complementary incomplete Gamma function, and the last equality follows from the well-known asymptotics for this function (see e.g.~\cite{olver97}, p.~66). On the other hand,
\[
\ell_p\alpha_pe^{-(\alpha_pJ)^2}\leq \ell_p\alpha_p e^{-(\log^{1/2}\ell_p\alpha_p-\alpha_p)^2}\leq e^{2\alpha_p\log^{1/2}(\ell_p\alpha_p)}\rightarrow 1.
\]
Therefore, the integral is bounded from above by a quantity of order lower than $\log^{3/2}(\ell_p\alpha_p)/\alpha_p$. To summarize,
\begin{equation}
 \label{log upper bound}
\sum_{j=0}^\infty\log\left(1+\ell_p \alpha_pe^{-(\alpha_p j)^2}\right)\leq \frac{2}{3}\frac{\log^{3/2}(\ell_p\alpha_p)}{\alpha_p}(1+o(1)).   
\end{equation}
Combining \eqref{log lower bound} and \eqref{log upper bound} yields \eqref{log asymptotics}.

Next, let now 
\[
J=\left\lfloor\frac{\log^{1/2}(a\ell_p\alpha_p)}{\alpha_p}\right\rfloor,
\]
where $a>0$ is an arbitrarily small constant. Such a $J$ is well defined for all sufficiently large $p$. 
We have $(\ell_p\alpha_p)^{-1}e^{(\alpha_p J)^2}\leq a$, and therefore
\begin{equation}
\label{lower bound der log}
\sum_{j=0}^{\infty}\frac{1}{(\ell_p \alpha_p)^{-1}e^{(\alpha_pj)^2}+1}\geq \frac{1+J}{1+a}\geq \frac{\log^{1/2} (a\ell_p \alpha_p)}{\alpha_p (1+a)}.
\end{equation}
On the other hand,
\begin{eqnarray*}
    \sum_{j=0}^{\infty}\frac{1}{(\ell_p \alpha_p)^{-1}e^{(\alpha_pj)^2}+1}&\leq&1+J+ \ell_p\alpha_p\sum_{j=J+1}^{\infty} e^{-(\alpha_pj)^2}\\
    &\leq& 1+J+\ell_p\alpha_p\int_{J}^{\infty} e^{-(\alpha_px)^2}\mathrm{d}x\\
    &=&1+J+\ell_p\sqrt{\pi}(1-\Phi(\sqrt{2}\alpha_pJ))\\
    &\leq&1+J+\frac{\ell_p}{2\alpha_pJ}e^{-\alpha_p^2J^2}.
\end{eqnarray*}
But
\begin{eqnarray*}
\frac{\ell_p}{J}e^{-\alpha_p^2 J^2}&\leq&(1+o(1))\frac{\ell_p\alpha_pe^{-\alpha_p^2\left(\log^{1/2}(a\ell_p\alpha_p)/\alpha_p-1\right)^2}}{\log^{1/2}(a\ell_p\alpha_p)}\\
&\leq&(1+o(1))\frac{e^{2\alpha_p\log^{1/2}(a\ell_p\alpha_p)}}{a\log^{1/2}(a\ell_p\alpha_p)}=\frac{1+o(1)}{a\log^{1/2}(a\ell_p\alpha_p)}=o(1).
\end{eqnarray*}
Therefore, 
\begin{equation}
    \label{upper bound der log}
\sum_{j=0}^{\infty}\frac{1}{(\ell_p \alpha_p)^{-1}e^{(\alpha_pj)^2}+1}\leq \frac{\log^{1/2}(a\ell_p\alpha_p)}{\alpha_p}(1+o(1)).
\end{equation}
Inequalities \eqref{lower bound der log}, \eqref{upper bound der log} and the fact that $a$ is an arbitrarily small positive number yield \eqref{log der asymptotics}.

\subsection*{Proof of Lemma \ref{lem: Stechkin 1}}

    Our proof of Lemma \ref{lem: Stechkin 1} closely follows the proof of a lemma on page 1511 of \cite{stechkin48}. We make only a few minor changes. Consider trigonometric polynomial $\Psi_k(\eta)=T_k(\bar{\omega}+\eta)-M\cos^2(k\eta/2)$. It has a double root at $\eta=0$, and, if  $T_k(\bar{\omega}+\eta)<M\cos^2\frac{k\eta}{2}$ for some $0<\eta<\pi/k$, it also has at least one additional root (the roots are counted according to their multiplicity) on each of the segments $\left[j\pi/k,(j+1)\pi/k\right]$ with $j=0,\dots,2k-2$. This would imply that $\Psi_k(\eta)$ has $2+2k-1=2k+1$ roots, which is impossible because its degree is $k$. The possibility that $T_k(\bar{\omega}+\eta)<M\cos^2\frac{k\eta}{2}$ for some $-\pi/k<\eta<0$ is eliminated similarly.

\end{appendix}
\bibliography{ddbib.bib}       


\end{document}